\documentclass[reqno]{amsart}

\usepackage{amscd,amsthm,amsmath,amssymb,mathtools}

\usepackage{upgreek,bm}

\usepackage{verbatim}
\usepackage{cite}
\usepackage{xspace}
\usepackage{xcolor}
\usepackage{url}
\usepackage{enumerate,enumitem}
\usepackage{graphicx}
 \usepackage{algorithm}
 \usepackage{algorithmic}

\usepackage{epstopdf,epsfig,subfigure}
\usepackage{curve2e}
\usepackage{setspace}
\usepackage{booktabs}

\usepackage[
    bookmarks=true,         
    unicode=false,          
    pdftoolbar=true,        
    pdfmenubar=true,        
    pdffitwindow=false,     
    pdfstartview={FitH},    
    pdftitle={Stability of RHC with positive semidefinite state costs},    
    pdfauthor={Author},     
    pdfsubject={Subject},   
    pdfcreator={Creator},   
    pdfproducer={Producer}, 
    pdfkeywords={keyword1, key2, key3}, 
    pdfnewwindow=true,      
    colorlinks=true,       
    linkcolor=red,          
    citecolor=blue,        
    filecolor=magenta,      
    urlcolor=cyan,           
hypertexnames=false
]{hyperref}

\usepackage{setspace}
\usepackage{geometry}
  \newcommand{\todoautc}[3][]{%
     \ifthenelse{\equal{#1}{}}{\todo[size=\scriptsize]{{\bf#2} #3}{}}{\todo[color=#1,size=\scriptsize]{{\bf#2} #3}{}}%
 }

 \newcommand{\noteC}[1]{\todoautc[cyan!50]{}{#1}}

\newcommand{\addref}[1][]{{\red\tt addref}\noteC{....}}

\theoremstyle{plain}

\newtheorem{theorem}{Theorem}[section]
\newtheorem{proposition}[theorem]{Proposition}
\newtheorem{lemma}[theorem]{Lemma}

\newtheorem{assumption}[theorem]{Assumption}

\theoremstyle{definition}
\newtheorem{definition}[theorem]{Definition}
\newtheorem{remark}[theorem]{Remark}

\numberwithin{equation}{section}

\newcommand{\linspan}{\mathop{\rm span}\nolimits}

\newcommand{\rest}{\left.\kern-2\nulldelimiterspace\right|_}
\newcommand{\norm}[2]{\left|#1\right|_{#2}}
\newcommand{\dnorm}[2]{\left\|#1\right\|_{#2}}

\newcommand{\argmin}{\mathop{\rm argmin}\nolimits}

\newcommand{\Id}{{\mathbf1}}
\newcommand{\indf}{1}

\newcommand{\p}{\partial}

\newcommand{\clA}{{\mathcal A}}
\newcommand{\clB}{{\mathcal B}}

\newcommand{\clE}{{\mathcal E}}

\newcommand{\clI}{{\mathcal I}}
\newcommand{\clJ}{{\mathcal J}}

\newcommand{\clL}{{\mathcal L}}

\newcommand{\clN}{{\mathcal N}}

\newcommand{\clP}{{\mathcal P}}

\newcommand{\clW}{{\mathcal W}}

\newcommand{\clZ}{{\mathcal Z}}

\newcommand{\bbN}{{\mathbb N}}

\newcommand{\bbR}{{\mathbb R}}

\newcommand{\bbU}{{\mathbb U}}

\newcommand{\bfM}{{\mathbf M}}

\newcommand{\bfS}{{\mathbf S}}

\newcommand{\fkI}{{\mathfrak I}}

\newcommand{\fkK}{{\mathfrak K}}

\newcommand{\fkT}{{\mathfrak T}}

\newcommand{\fkY}{{\mathfrak Y}}

\newcommand{\rmD}{{\mathrm D}}

\newcommand{\bfn}{{\mathbf n}}

\newcommand{\rmd}{{\mathrm d}}
\newcommand{\rme}{{\mathrm e}}
\newcommand{\rmf}{{\mathrm f}}

\definecolor{DarkBlue}{rgb}{0,0.08,0.45}
\definecolor{DarkRed}{rgb}{.65,0,0}
\definecolor{applegreen}{rgb}{0.55, 0.71, 0.0}

\newcounter{mymac@matlab}
\newcommand{\matlab}{MATLAB%
   \ifnum\value{mymac@matlab}<1%
   \textregistered%
   \setcounter{mymac@matlab}{1}%
   \fi%
  }

\providecommand{\argmin}{\operatorname*{argmin}}

\newcommand{\red}{ \color{red} }

 \usepackage[mathlines,pagewise]{lineno}

 \let\oldequation\equation
 \let\oldendequation\endequation

 \renewenvironment{equation}
   {\linenomath\oldequation}
   {\oldendequation\linenomath}

 \let\oldalign\align
 \let\oldendalign\endalign

 \renewenvironment{align}
   {\linenomath\oldalign}
   {\oldendalign\linenomath}

\begin{document}
\title{Receding horizon control with sensor-induced positive semidefinite state penalization}
\author{Behzad Azmi$^{\tt1}$}
\author{Sergio S.  Rodrigues $^{\tt2}$}

 \thanks{
\vspace{-1em}\newline\noindent
{\sc MSC2020}: 93C05, 93C10, 93C20, 93D20\quad\hfill
\newline\noindent
$^{\tt1}$ Department of Mathematics and Statistics, University of Konstanz, D-78457 Konstanz, Germany \newline\noindent
$^{\tt2}$ Department of Mathematics,
FCT, NOVA
University of Lisbon,
2829-516 Caparica,
Portugal.
\newline\noindent
  {\sc Emails}:
 {\small\tt behzad.azmi@uni-konstanz.de,\quad ssi.rodrigues@fct.unl.pt}.
  \newline\noindent
}

\begin{abstract}

This paper develops and analyzes a receding horizon control (RHC) framework for the stabilization of infinite-dimensional nonautonomous continuous-time systems with an infinite-horizon control objective. The considered stage costs penalize only a finite-dimensional component of the state and are therefore positive semidefinite with respect to the full state. Without imposing terminal costs or terminal constraints, we establish exponential stability and derive suboptimality estimates under suitable stabilizability and detectability-type assumptions, together with an online selection of the sampling times. In the linear case, we extend the framework to output-feedback RHC by combining it with a Luenberger-type state estimator based
on finitely many sensors, and establish closed-loop stability using an augmented Lyapunov functional. The abstract state-feedback results are verified for a nonlinear Schl\"ogl parabolic equation around an unstable steady state and, in the output-feedback setting, for nonautonomous linear parabolic equations. The control and observation realizations involve finite-dimensional actuation and sensing, including localized boundary controls and spatially localized state-penalization operators. Numerical experiments support the theoretical findings and illustrate the performance of the proposed RHC schemes.
\end{abstract}

\maketitle

\pagestyle{myheadings} \thispagestyle{plain} \markboth{\sc B. Azmi and S. S.
Rodrigues}{\sc stability of RHC with positive semidefinite state penalization operators}

\section{Introduction}

Many processes in science and engineering are modeled by parabolic partial
differential equations, including heat transfer, diffusion, reaction--diffusion
phenomena, and phase-transition dynamics. From a control perspective, these
systems are challenging because their state spaces are infinite-dimensional,
the control action is often available only through finitely many localized
actuators, and the state information may be restricted to finitely many
measurements. These features make parabolic control problems a natural setting
for feedback strategies that can handle large-scale dynamics, localized
actuation, and partial observation.

We are concerned with controlled systems of the form
\begin{equation}
\label{cs}
    \dot y + Ay + \mathcal{N}(\cdot,y) = Bu,
    \qquad
    y(0)=y_0\in Y,
    \qquad
    t\in(0,\infty) .
\end{equation}
Here $y(t)\in Y$ denotes the state, evolving in a Hilbert space $Y$, and
$u(t)\in\bbR^m$ is the control input. The operator
$A\colon \rmD(A)\subset H\to H$ is linear and time-independent, while
$\mathcal{N}\colon [0,\infty)\times X\to Z$ is a possibly
time-dependent and possibly nonlinear operator, where $X\subset Y\subset Z$
are suitable Hilbert spaces. We assume that $\mathcal{N}(\cdot,0)=0$, and
$B\in\mathcal{L}(\bbR^m,Z)$ denotes the control operator. The
finite-dimensionality of the input reflects the fact that, in many
applications, only a finite number of actuators is available.

The control objective is to design a feedback mechanism which stabilizes the
origin. More precisely, we seek controls
$u\in L^2(\bbR_+,\bbR^m)$ such that the corresponding solution satisfies an
exponential decay estimate of the form
\[
    \norm{y^{u}(t)}{Y}
    \leq
    C e^{-\alpha t}\norm{y_0}{Y},
    \qquad t\geq 0,
\]
for some constants $C\geq 1$ and $\alpha>0$.

Receding horizon control (RHC) provides a flexible feedback strategy for this
purpose. At each sampling time, a finite-horizon optimal control problem is
solved, only the first part of the computed control is applied, and the
procedure is repeated. The resulting closed-loop performance is then measured
with respect to the corresponding infinite-horizon cost. In this paper, we
focus on stage costs of the form
\[
    \ell(y(t),u(t))
    =
    |Qy(t)|_{\bbR^s}^2
    +
    |u(t)|_{\bbR^m}^2,
\]
where $Q:Y\to\bbR^s$ is a finite-dimensional state-penalization operator.
Since $Q$ has a finite-dimensional range, this stage cost is only positive
semidefinite with respect to the state and does not define a norm on the full
state space. Consequently, stability cannot be inferred from cost alone
and must be recovered through suitable stabilizability and detectability
properties of the dynamics and the operator $Q$.

The use of such finite-dimensional state penalties is motivated by robustness,
performance, and computational considerations. Full-state penalties may
increase the sensitivity to model errors, discretization effects, and measurement
noise, especially in poorly resolved parts of the state. Moreover, many
control objectives are naturally output-oriented, for instance, involving
averaged quantities, localized measurements, fluxes, interface positions, or
trace functionals. In such cases, $Q$ can be interpreted as an observation or
output operator, and the term $|Qy|_{\bbR^s}^2$ focuses the control effort on
the relevant performance variables. Finally, finite-dimensional state
penalization is attractive in online RHC, since the state contribution to the
cost is evaluated in a low-dimensional space rather than on the full spatial
discretization. Nevertheless, stability can be recovered under appropriate
detectability assumptions.

The continuous-time character of the problem introduces an additional
difficulty. In contrast to discrete-time RHC frameworks, the detectability
mechanism used here is formulated over continuous-time sampling intervals.
This is important for sensor-induced positive semidefinite stage costs, where
the relevant observability or detectability constants may depend on the
sampling length, the parameters of the equation, and the particular choice of
the operator $Q$. Hence, these constants cannot, in general, be fixed uniformly
for arbitrary a priori sampling times. For this reason, the sampling times in
the RHC scheme considered here are not fixed beforehand, but are selected
online through a synchronization procedure.

\subsection{Literature review}

Receding horizon control (RHC), also known as model predictive control, is a
well-established framework for addressing infinite-horizon optimal control
problems through the repeated solution of finite-horizon problems; see, for
example, \cite{GruenePannek17,RawlingsMayneDiehl19}. Its stability analysis
has a long history in finite-dimensional systems, both in continuous and
discrete time; see, among many others,
\cite{MayneMichalska90,MayneRawlRaoScok00,GruenePannek09,GrueneRantzer08,Jad01}
and the references therein. While classical Riccati-based feedback
constructions are available for important linear-quadratic problems, RHC is
particularly useful in settings involving constraints, nonlinearities,
non-autonomous dynamics, or large-scale systems, where direct
infinite-horizon synthesis is often difficult.

Several recent works have developed stability and performance estimates for
RHC schemes in the setting of finite-dimensional discrete-time control
systems. In particular, positive semidefinite stage costs can be treated under
suitable detectability conditions; see, for instance,
\cite{Grimm05model,Koehler23stability,KohMulAll22}. These results are
closely related in spirit to the present work, since they also address the
difficulty caused by stage costs that do not directly measure the full state.
However, their detectability assumptions are tailored to discrete-time
dynamics, where the relevant estimate is imposed between two consecutive time
instants. In the continuous-time parabolic setting considered here, the
corresponding detectability mechanism is a finite-time property; its constants
may depend on the sampling interval, the system parameters, and the chosen
observation operator.

RHC has also been studied in infinite-dimensional and PDE-related settings.
Ito and Kunisch developed RHC approaches for systems described in function
spaces, including formulations with terminal ingredients and incomplete
observations; see \cite{ItoKun02,ItoKun06}. Stabilizing RHC
schemes under detectability-type assumptions and terminal costs were further
investigated in \cite{KunPfeiffer20}. Terminal-cost-free RHC schemes for
infinite-dimensional systems were studied, for instance, in
\cite{AzmiKunisch19,azmi22}. The present work contributes to this line of research
by treating abstract semilinear continuous-time systems with sensor-induced
positive semidefinite stage costs. A distinctive feature of our framework is
that the sampling times are selected online through a synchronization
mechanism, rather than being fixed a priori.

The output-feedback part of the paper is related to observer-based RHC. For
finite-dimensional systems, the combination of RHC with state estimators has
been studied in several settings; see, for instance,
\cite{MayneRavicFindAllg09,FinImsAllgFos03,KohMulAll22}. These works address partial-state information, estimation
errors, and robustness issues by combining predictive control with suitable
observer or estimator dynamics. In our previous work~\cite{Azmi25output}, an
output-feedback RHC (ORHC) algorithm with a Luenberger-type state estimator was
introduced for parabolic systems. In the present paper, we use the same
algorithmic structure in the linear case, but develop a different stability
analysis. More precisely, the proof is based on an augmented Lyapunov
functional that couples the RHC closed-loop dynamics with the observer error
dynamics. This provides a direct stability argument for the ORHC scheme in the present infinite-dimensional setting.

\subsection{Contributions}

The main contributions of this paper can be summarized as follows.

\begin{itemize}
    \item We develop a receding horizon control (RHC) framework for a class of
    abstract semilinear evolution equations in continuous time. The analysis
    allows for stage costs that are only positive semidefinite, and hence not
    positive definite with respect to the full state. This is particularly
    relevant for infinite-dimensional systems where the cost may only measure
    finitely many state components or outputs.

    \item For this class of systems, we prove stability and suboptimality
    properties of the proposed RHC scheme. More precisely, under suitable
    stabilizability- and detectability-type assumptions, we establish
    exponential stability of the closed-loop system and derive a corresponding
    suboptimality estimate. A key feature of the analysis is that the sampling
    times are not fixed a priori, but are selected online in accordance with
    the synchronization mechanism of the algorithm.

    \item In the linear case, we revisit an ORHC framework that
    was introduced in our previous work~\cite{Azmi25output}. While the
    algorithmic structure is retained, the stability analysis developed here
    is new. In particular, we combine the RHC scheme with a Luenberger-type
    state estimator based on finitely many sensors and prove stability of the
    resulting observer-based closed loop by means of an augmented Lyapunov
    functional. This provides a different and more direct analysis of the
    ORHC scheme.

    \item We apply the abstract framework to nonlinear Schl\"ogl equations
    with finite-di\-men\-sion\-al internal or dynamic boundary control, and to
    nonautonomous linear parabolic equations in the output-feedback setting.
    These applications involve localized actuators and sensors, with positive
    semidefinite state costs generated by finitely many spatial measurements.
 
    \item Finally, we provide comparative numerical simulations for the
    full-state and ORHC schemes. The simulations support the
    theoretical stability and suboptimality results and illustrate the
    effectiveness of the proposed online sampling strategy for systems with
    localized actuation and sensing.
\end{itemize}

\subsection{Contents}

The paper is organized as follows. Section~\ref{S:stabil_rhc} introduces
the abstract continuous-time RHC framework and shows the main results of stabilizability and suboptimality. Section~\ref{sec:output_feedback} treats the ORHC scheme in the linear case.
Section~\ref{sec:parabolic_verification} verifies the abstract assumptions
for classes of parabolic equations, in several settings including localized finite-dimensional boundary actuation and sensing. Numerical examples are presented in
Section~\ref{sec:numerics} for the case of boundary actuation and internal sensing. Concluding remarks are given in Section~\ref{sec:conclusion}. Some technical results used in the optimal-control analysis are provided in Appendix~\ref{app:schl-wellposedness}.

\subsection*{Notation}
We write~$\bbR$ and~$\bbN$ for the sets of real numbers and nonnegative
integers, respectively, and set
$\bbR_+\coloneqq(0,+\infty)$,
$\overline{\bbR}_+\coloneqq[0,+\infty)$, and
$\bbN_+\coloneqq\bbN\setminus\{0\}$.

Given Hilbert spaces~$ X$ and~$Y$, the space of continuous linear mappings from~$ X$ into~$Y$ is denoted by~$\clL(X,Y)$. If $X=Y$, we write $\clL(X)\coloneqq\clL(X,X)$. The continuous dual of~$X$ is denoted by
$X'\coloneqq\clL(X,\bbR)$. The adjoint of an operator
$L\in\clL(X,Y)$ is denoted by~$L^*\in\clL(Y',X')$.
If the embedding $X\subseteq Y$ is continuous, we write
$X\xhookrightarrow{}Y$. We write $X\xhookrightarrow{\rm d}Y$ and
$X\xhookrightarrow{\rm c}Y$ if the embedding is, in addition, dense and
compact, respectively. For $z_0\in Y$ and $r>0$, we denote by
$\clB_r^Y(z_0)\coloneqq\{z\in Y\mid \norm{z-z_0}{Y}<r\}$ the open ball in
$Y$ centered at~$z_0$ with radius~$r$.

For an open interval~$I\subset\bbR$, we denote by $W(I;X,Y)\coloneqq\{f\in L^2(I;X)\mid \dot f\in L^2(I;Y)\}$ the function space endowed with the norm
$\norm{f}{W(I;X,Y)}
\coloneqq
\norm{(f,\dot f)}{L^2(I;X)\times L^2(I;Y)}$.

Finally, for $\overline t_0\in\overline{\bbR}_+$ and $T\in\bbR_+$, we set
$\bbR_{\overline t_0}\coloneqq(\overline t_0,\infty)$ and
$\clI_{\overline t_0,T}\coloneqq(\overline t_0,\overline t_0+T)$, with
closures denoted by
$\overline{\bbR}_{\overline t_0}\coloneqq[\overline t_0,\infty)$ and
$\overline{\clI}_{\overline t_0,T}
\coloneqq[\overline t_0,\overline t_0+T]$.

\section{Stabilizing receding-horizon control}\label{S:stabil_rhc}

We now formulate the receding horizon control strategy for the controlled
system~\eqref{cs}. Throughout this section, we use the notation introduced
above and assume that, for every admissible control and every initial pair
$(\overline t_0,\overline y_0)\in\overline{\bbR}_+\times Y$, the corresponding solution is well defined on the time interval under
consideration. For $T\in\bbR_+\cup\{+\infty\}$,
$(\overline t_0,\overline y_0)\in\overline{\bbR}_+\times Y$, and a control
$u$ defined on $\clI_{\overline t_0,T}$, we denote by
\[
    \fkY_T(\overline t_0,\overline y_0;u)(t)\coloneqq y(t)
\]
the solution of
\begin{equation}
\label{control_sys_b}
    \dot y + Ay+\mathcal{N}(\cdot,y)=Bu,
    \qquad
    y(\overline t_0)=\overline y_0,
    \qquad
    t\in\clI_{\overline t_0,T}
    \coloneqq
    (\overline t_0,\overline t_0+T).
\end{equation}
Thus, \eqref{control_sys_b} is the shifted finite-horizon version of
\eqref{cs}, starting at time $\overline t_0$ from the state
$\overline y_0$. For an admissible set $\bbU\subseteq\bbR^m$, we define
\[
    L^2_{\rm ad}(\bbR_+,\bbR^m)
    \coloneqq
    \{v\in L^2(\bbR_+,\bbR^m)\mid v(t)\in\bbU
    \ \text{for a.e. } t\in\bbR_+\}.
\]
Analogously, for a finite time interval $\clI_{\overline t_0,T}$, we set
\[
    L^2_{\rm ad}(\clI_{\overline t_0,T},\bbR^m)
    \coloneqq
    \{v\in L^2(\clI_{\overline t_0,T},\bbR^m)
    \mid v(t)\in\bbU
    \ \text{for a.e. } t\in\clI_{\overline t_0,T}\}.
\]
For example, $\bbU$ may be chosen as
\[
    \bbU=\{v\in\bbR^m\mid \dnorm{v}{}\leq C_u\},
\]
for some norm $\dnorm{\cdot}{}$ on $\bbR^m$ and some
$0<C_u\leq\infty$. The case $C_u<\infty$ corresponds to pointwise magnitude
constraints on the input, while $C_u=\infty$ recovers the unconstrained case. The finite-horizon performance index is then defined by
\begin{subequations}\label{Obj}
\begin{align}
    J_T(u;\overline t_0,\overline y_0)
    &\coloneqq
    \int_{\overline t_0}^{\overline t_0+T}
    \ell\left(
        \fkY_T(\overline t_0,\overline y_0;u)(t),u(t)
    \right)\,\rmd t,
    \\
    \ell(z,v)
    &\coloneqq
    \norm{Qz}{\bbR^s}^2+\norm{v}{\bbR^m}^2 .
\end{align}
\end{subequations}
Here $Q\in\clL(Y,\bbR^s)$ is the finite-dimensional state-penalization
operator introduced above. Since $Q$ has finite-dimensional range, the state
contribution in $\ell$ is only positive semidefinite with respect to the full
state. For a given initial state $y_0\in Y$, the corresponding infinite-time-horizon
problem is formally given by
\begin{align}
\label{tag:opinf}\stepcounter{equation}
 \tag*{{\rm (\theequation:\text{$P_{\infty}^{0,y_0}$})}}
    \min_{u\in L^2_{\rm ad}(\bbR_+,\bbR^m)}
    J_\infty(u;0,y_0).
\end{align}
The receding horizon strategy approximates the feedback effect of this
infinite-horizon problem through the concatenation of controls obtained from
a sequence of finite-time-horizon optimal control problems. More precisely,
for a current time-state pair $(\overline t_0,\overline y_0)$ and a
prediction horizon $T>0$, we consider
\begin{align}
\label{tag:opfin}\stepcounter{equation}
 \tag*{{\rm (\theequation:\text{$P_T^{\overline t_0,\overline y_0}$})}}
    \min_{u\in L^2_{\rm ad}(\clI_{\overline t_0,T},\bbR^m)}
    J_T(u;\overline t_0,\overline y_0).
\end{align}
Whenever it exists, an optimal solution of~\ref{tag:opfin} is denoted by
$u_T^{\overline t_0,\overline y_0,*}$. For later use, we introduce the notation
\[
    \ell_T^*(t;\overline t_0,\overline y_0)
    \coloneqq
    \left(
        \norm{
            Q\fkY_T(
                \overline t_0,\overline y_0;
                u_T^{\overline t_0,\overline y_0,*}
            )(t)
        }{\bbR^s}^2
        +
        \norm{
            u_T^{\overline t_0,\overline y_0,*}(t)
        }{\bbR^m}^2
    \right),
    \qquad
    t\in\clI_{\overline t_0,T},
\]
for the integrand of the optimal finite-horizon cost. Hence,
\[
    J_T(u_T^{\overline t_0,\overline y_0,*};
        \overline t_0,\overline y_0)
    =
    \int_{\overline t_0}^{\overline t_0+T}
    \ell_T^*(t;\overline t_0,\overline y_0)\,\rmd t .
\]
We use the same normalization for the infinite-horizon cost. Omitting the
common prefactor $1/2$ does not change the minimizers; in particular,
$\norm{u}{L^2(\clI_{\overline t_0,T};\bbR^m)}^2
\le J_T(u;\overline t_0,\overline y_0)$.
In Algorithm~\ref{alg:RH}, we recall the receding horizon strategy used in this work.
\begin{algorithm}[htbp]
\caption{Receding horizon control,   RHC($\tau,T$)}\label{alg:RH}
\begin{algorithmic}[1]
\REQUIRE{lower bound $\tau$ for sampling time; prediction horizon $T\geq \tau$; initial state $y_0\in Y$.}
\ENSURE{a control input~$u_{\rm rh}\in L^2_{\rm ad}( \bbR_+ , \bbR^m)$; concatenation times~$(t_{k})_{k\in\bbN}$.}
\STATE Set~$(\overline t_0,\overline{y}_0):=(0, y_0)$, $t_0=0$,  and  $y_{\rm rh}(0) =y_0 $;
\STATE\label{alg:Tinf} Set~$T_{\infty}=\infty$; ${k}=0$;
\WHILE{$\overline t_0<T_{\infty}$}
\STATE\label{alg:solve} ``Solve'' the  FTH problem~\ref{tag:opfin} to find~$u^{\overline t_0,\overline{y}_0,*}_T  \in L^2(\clI_{\overline{t}_0,T},U)$;
\STATE\label{alg:delta} Set  $\delta :=\max \argmin_{s \in [\tau, T]} \norm{\fkY_T(\overline t_0,\overline y_0;u^{\overline t_0,\overline{y}_0,*}_T)( \overline{t}_0+s)}{Y}^2$;
\STATE For all $t \in [\overline t_0, \overline t_0+\delta)$,  set $u_{\rm rh}(t)= u^{\overline t_0,\overline{y}_0,*}_T(t)$ and~$y_{\rm rh}(t)=\fkY_T(\overline t_0,\overline{y}_0;u^{\overline t_0,\overline{y}_0,*}_T)(t)$;
\STATE\label{alg:RH-update-state} Update: $(\overline t_0,\overline{y}_0)  \leftarrow (\overline t_0+\delta,  \fkY_T(\overline t_0,\overline{y}_0;u^{\overline t_0,\overline{y}_0,*}_T)(\overline t_0+\delta))$;
\STATE Set ${k}\gets {k}+1$;\quad $t_{k}\gets \overline t_0$;
\ENDWHILE
\end{algorithmic}
\end{algorithm}
At each iteration, a finite-horizon problem of the form \ref{tag:opfin} is solved, the computed optimal control is applied only on the current sampling interval, and the procedure is repeated from the newly reached state. The notation $T_\infty$ is introduced in Algorithm~\ref{alg:RH} for later
convenience. In the numerical simulations, we will use a finite value
$T_\infty<\infty$ as a stopping criterion for the RHC procedure.

\subsection{Stabilizability and suboptimality}\label{sS:rhc-stabsubopt}
In this section,  we show that the receding horizon control provided by Algorithm \ref{alg:RH} is stabilizing. Furthermore, we investigate the suboptimality of its associated cost compared to the optimal cost of the infinite-horizon problem. To do this, we need to define the so-called value functions, giving us the optimal (i.e., minimal) cost of the involved optimal control problems.

We will follow a variation of the strategy developed in \cite{AzmiKunisch19} for positive definite~$Q$; see the assumption in~\cite[Eq.~(2.3)]{AzmiKunisch19}. Since the operator  $Q$ in \eqref{Obj} is not positive definite, we will need some new arguments. In particular, we need to show the dissipative inequality \eqref{e20s} by relying on requirement~\ref{asm:rh-estimYT} (in particular,~\eqref{e8}) within Assumption~\ref{asm:rh}, which is related to the detectability condition for $Q$.
For any $y_0 \in Y$,  ITH value function $V_{\infty}$ is defined by
\begin{subequations}\label{eq.valuefunct}
\begin{equation}
 V_{\infty}(y_0)\coloneqq \inf_{u \in L^2_{\rm ad}(\bbR_+ , \bbU)} J_{\infty}(u;  0, y_0).
\end{equation}
Similarly, for every initial vector  $(\overline{t}_0,  \overline{y}_0) \in \overline\bbR_{+} \times Y$ and every prediction horizon $T>0$,  the FTH value function $V_{T}$ is defined by
\begin{equation}
V_{T}(\overline{t}_0,\overline{y}_0) \coloneqq \mathop{\inf}\limits_{u\in L^2_{\rm ad}(\clI_{\overline{t}_0,T}, \bbU)} J_{T}(u;\overline{t}_0,\overline{y}_0).
\end{equation}
\end{subequations}
Under requirement~\ref{asm:rh-existopt} below, this infimum is attained for the horizons and initial data covered by that assumption.

To demonstrate the exponential stability and suboptimality of the RHC input obtained by Algorithm \ref{alg:RH}, we need some properties as in the following assumption.
In particular, we need to assume/guarantee that the involved optimal control problems are well-posed.

\begin{assumption}\label{asm:rh}
    With~$\tau>0$ being the minimal sampling time~$\tau$ in Algorithm~\ref{alg:RH}, there exists a radius $r>0$ such that:
    \begin{enumerate}[noitemsep,topsep=5pt,parsep=5pt,partopsep=0pt,leftmargin=0em]%
\renewcommand{\theenumi}{{\sf RH\arabic{enumi}}} 
 \renewcommand{\labelenumi}{}
 \item\theenumi{\rm:}\label{asm:rh-existopt} For every  $(T,\overline{t}_0 ,\overline{y}_0) \in \overline\bbR_{\tau} \times\overline\bbR_{+} \times\clB^Y_{r}(0)$, the optimal control problem \ref{tag:opfin}  admits a solution.
\item\theenumi{\rm:}\label{asm:rh-bddVT} There exists a continuous,  nondecreasing,  bounded function $\gamma: \overline{\mathbb{R}}_\tau \to \overline{\mathbb{R}}_+$ such that   
\begin{equation}
\label{e7}
 V_T( \overline{t}_0, \overline{y}_0)  \leq  \gamma(T)\norm{ \overline{y}_0}{Y}^2 \quad \text{ for every } (\overline{t}_0, \overline{y}_0)\in   \overline\bbR_{+} \times\clB^Y_{r}(0).
\end{equation}
 \item\theenumi{\rm:}\label{asm:rh-estimYT} For each~$T\ge\tau$, there are constants~$\overline{C}(T)\ge1$ and~${C}_{\rm det}(T)\ge1$ with~$\sup_{T\ge\tau}{C}_{\rm det}(T)<\infty$ such that for every  $(\overline{t}_0 ,\overline{y}_0,u) \in \overline\bbR_{+} \times\clB^Y_{r}(0)\times L^2_{\rm ad}(\clI_{\overline{t}_0,T},\bbR^m)$, satisfying 
\begin{equation}
\label{contC}
 \norm{u}{L^2(\clI_{\overline{t}_0,T},\bbR^m)}\leq \sqrt{\gamma(T)}\, \norm{\overline{y}_0}{Y},
\end{equation}
 there exists one, and only one,  solution for~\eqref{control_sys_b}. Moreover,
\begin{align}
\label{e27}
&\norm{\fkY_T(\overline t_0,\overline y_0; u)}{C(\overline{\clI}_{\overline{t}_0,T}, Y)}^2 \leq  \overline{C}(T) \left(  \norm{\overline{y}_0}{Y}^2+ \norm{u}{L^2(\clI_{\overline{t}_0,T},\bbR^m)}^2 \right);\\
&\norm{ \fkY_T(\overline t_0,\overline y_0; u)}{L^2(\clI_{\overline{t}_0,T} ;Y)}^2  \leq C_{\rm det}(T) \left(  \norm{\overline{y}_0}{Y}^2 +  J_{T}(u; \overline t_0,  \overline{y}_0)  \right).\label{e8}
\end{align}  
     \end{enumerate}
\end{assumption}

\begin{remark}
The stability and performance analysis in
\cite{Koehler23stability} for discrete-time systems relies, among other
hypotheses, on a detectability condition; see, for instance,
\cite[Assumption~3]{Koehler23stability},
\cite[SA3--SA4]{Grimm05model}, and
\cite[Assumptions~3--4]{KohMulAll22}. In common notation, this condition
requires the existence of a nonnegative storage function $W$ and constants
$\gamma_o \geq 0$ and $\varepsilon_o>0$ such that
\[
    W(x) \leq \gamma_o \sigma(x),
    \qquad
    W(f(x,u))-W(x)
    \leq
    -\varepsilon_o \sigma(x)+\ell(x,u),
\]
for all admissible states $x$ and inputs $u$. This is a one-step condition,
naturally formulated for discrete-time dynamics.
For a continuous-time system, the corresponding formulation is more naturally
given over a sampling interval. Even in the autonomous case, if
$\varphi(t;x,u)$ denotes the solution starting from $x$ under the control
$u$, one is led to estimates of the form
\begin{equation}\label{arb-small-deltaQ}
    W(\varphi(\delta;x,u))-W(x)
    \leq
    -\varepsilon_\delta \sigma(x)
    +
    \int_0^\delta \ell(\varphi(t;x,u),u(t))\,dt.
\end{equation}
Here,
the constant $\varepsilon_\delta$ generally depends on the sampling
length $\delta$. In particular, when the stage cost is only positive
semidefinite, the information not directly measured by the cost has to be
recovered through the dynamics over time. For parabolic systems, this is
typically a finite-time observability or detectability mechanism. Hence, even
in the autonomous case, it is not clear whether~\eqref{arb-small-deltaQ} will be satisfied for arbitrarily small~$\delta$. Thus, the sampling interval may need 
to be sufficiently large, depending on
the equation, its parameters, and the chosen observation or stage cost.
The situation is even more delicate for non-autonomous systems, since the flow
depends explicitly on the initial time, say
$\varphi(\overline t_0+\delta;\overline t_0,x,u)$. The corresponding
detectability constants may then depend not only on the sampling length
$\delta$, but also on the initial time $\overline t_0$. Thus, unless a
suitable uniform-in-time detectability estimate is available, such constants
cannot generally be chosen a priori uniformly over all sampling intervals and
all initial times.
This issue is relevant in the present manuscript, where the stage cost is
induced by finitely many sensors and is therefore not positive definite with
respect to the full infinite-dimensional state. Consequently, the above
discrete-time detectability mechanism does not directly provide the uniform
estimates needed in our continuous-time setting for arbitrary a priori
prescribed sampling times. Under the present assumptions and proof strategy,
this is one reason why stability and suboptimality are proved with sampling
times selected online through the synchronization mechanism of the RHC scheme. Finally, note that in this manuscript we are interested in stabilization to zero, which would correspond, in the language in~\cite[above Def.~1]{Grimm05model}, to the stabilization of the set~$\clA=\{0\}$.
 
\end{remark}

Next, we define the augmented control Lyapunov function that will be used
throughout the stability analysis.

\begin{definition}[Augmented control Lyapunov function]
\label{def:CLF}
Let $\tau>0$ be fixed. We choose a function
$h$ satisfying
\begin{subequations}\label{clf}
\begin{align}
  &h\colon \mathbb{R}_\tau\to \mathbb{R}_+,\qquad \lim_{T\to\infty}h(T) = 0,\qquad\lim_{T\to\infty}  \frac{1}{h(T)(T-\tau)}=0.
  \intertext{For $T>\tau$, we define}
  &G_T(t,z)\coloneqq  V_T(t,z)+h(T) \norm{z}{Y}^2,\qquad(t,z)\in\overline\bbR_+\times \clB^Y_r(0).\label{clf-GT=VT+h}
\end{align}
\end{subequations}
We refer to $G_T$ as the augmented control Lyapunov function.
\end{definition}

\begin{lemma}
\label{lem:bdd_GT}
Let Assumption~\ref{asm:rh} hold for some $\tau>0$ and~$r>0$.
For every fixed $T>\tau$, there exists a radius $r_1=r_1(T)\le r$
such that, for every $(\overline{t}_0,\overline{y}_0)\in
\overline\bbR_+\times\clB^Y_{r_1(T)}(0)$,
with $G_T$ defined as in~\eqref{clf}, we have that
\begin{equation}
\label{e23}
\begin{split}
G_T( \overline{t}_0+\delta, \fkY_T(\overline t_0,\overline y_0;u^{\overline t_0,\overline{y}_0,*}_T)(\overline{t}_0+\delta)) \leq \tfrac{C_{\rm det}(T)(\gamma(T) + h(T) )}{T-\tau} \left( \norm{\overline{y}_0}{Y}^2+ V_T(\bar{t}_0, \bar{y}_0)\right), 
\end{split}
\end{equation}
where~$\delta\in[\tau,T]$ is a minimizer of
$s\mapsto\norm{\fkY_T(\overline t_0,\overline y_0;u^{\overline t_0,\overline{y}_0,*}_T)(\overline t_0+s)}{Y}$,
as in Algorithm~\ref{alg:RH}.
\end{lemma}
\begin{proof}
We set $u^*\coloneqq u_T^{\overline t_0,\overline y_0,*}$ and
$\overline\fkY\coloneqq\fkY_T(\overline t_0,\overline y_0;u^*)$.
Then, by the normalization in~\eqref{Obj} and~\ref{asm:rh-bddVT}, we obtain
\[
 \norm{u^*}{L^2(\clI_{\overline t_0,T};\bbR^m)}^2
 \le J_T(u^*;\overline t_0,\overline y_0)
 =V_T(\overline t_0,\overline y_0)
 \le\gamma(T)\norm{\overline y_0}{Y}^2.
\]
Consequently, $u^*$ satisfies~\eqref{contC}. Applying~\eqref{e27} with
horizon $T$, we obtain, for the sampling time $\delta\in[\tau,T]$
selected in Algorithm~\ref{alg:RH},
\[
 \norm{\overline\fkY(\overline t_0+\delta)}{Y}^2
 \le\norm{\overline\fkY}{C(\overline{\clI}_{\overline t_0,T};Y)}^2
 \le\overline C(T)\left(\norm{\overline y_0}{Y}^2
       +\norm{u^*}{L^2(\clI_{\overline t_0,T};\bbR^m)}^2\right)
 \le\overline C(T)(1+\gamma(T))\norm{\overline y_0}{Y}^2.
\]
Consequently, choosing
$r_1(T)\coloneqq r/\sqrt{\overline C(T)(1+\gamma(T))}\le r$,
we have $\overline\fkY(\overline t_0+\delta)\in\clB_r^Y(0)$
whenever $\overline y_0\in\clB_{r_1(T)}^Y(0)$.
Furthermore,  since  $\delta\in \argmin_{t \in [\tau, T]}  \norm{ \overline\fkY( \overline{t}_0+t)}{Y}^2$, we can infer that 
\begin{align}
 \norm{ \overline\fkY(\overline{t}_0+\delta)}{Y}^2 &\leq \tfrac{1}{T-\tau} \norm{ \overline\fkY}{L^2((\overline{t}_0+\tau,\overline{t}_0+T),Y)}^2  \leq \tfrac{1}{T-\tau}  \norm{\overline\fkY}{L^2((\overline{t}_0,\overline{t}_0+T),Y)}^2\notag\\
& \leq \tfrac{C_{\rm det}(T)}{T-\tau} \left( \norm{\overline{y}_0}{Y}^2+  \int_{ \overline{t}_0}^{\overline{t}_0+T} \ell^*_T(t;\overline{t}_0, \overline{y}_0)dt      \right)\notag
\end{align}
which, together with the definitions of $G_T$, in~\eqref{clf}, and of the value function~$V_T$, give us
\begin{align}
&G_T( \overline{t}_0+\delta, \overline\fkY(\overline{t}_0+\delta))=
V_T( \overline{t}_0+\delta, \overline\fkY(\overline{t}_0+\delta))+h(T) \norm{\overline\fkY(\overline{t}_0+\delta)}{Y}^2\notag
\\&\hspace{2em} \le(\gamma(T) + h(T)) \norm{\overline\fkY(\overline{t}_0+\delta)}{Y}^2\le \tfrac{C_{\rm det}(T)(\gamma(T) + h(T) )}{T-\tau} \left( \norm{\overline{y}_0}{Y}^2+  {\textstyle\int_{ \overline{t}_0}^{\overline{t}_0+T}} \ell^*_T(t;\overline{t}_0, \overline{y}_0)dt  \right)\notag\\
&\hspace{2em}  =  \tfrac{C_{\rm det}(T)(\gamma(T) + h(T) )}{T-\tau} \left( \norm{\overline{y}_0}{Y}^2+  V_T(\overline{t}_0, \overline{y}_0)\right),\notag
\end{align}
which completes the proof.
\end{proof}

\begin{lemma}
\label{lem:estGT}
Let Assumption~\ref{asm:rh} hold true and let $\tau>0$ and~$\alpha\in(0,1)$ be given. Then there exists $T^*>\tau$ such that for every $T\geq T^*$ and
$(\overline{t}_0,\overline{y}_0)\in\overline\bbR_+\times\clB^Y_{r_1(T)}(0)$,
with $r_1(T)$ defined in Lemma~\ref{lem:bdd_GT},
we have for $G_T(t,y)$ defined as in~\eqref{clf} that
\begin{equation}
\label{e20s}
G_T(\overline{t}_0+\delta ,\fkY_T(\overline t_0,\overline y_0;u^{\overline t_0,\overline{y}_0,*}_T)(\overline{t}_0+\delta)) - G_T(\overline{t}_0, \overline{y}_0)   \leq    -\alpha \int_{\overline{t}_0}^{\overline{t}_0+\delta} \ell^*_T(t;\overline{t}_0, \overline{y}_0)dt 
\end{equation}
and   
\begin{equation}
\label{e20ss}
G_T(\overline{t}_0+\delta ,\fkY_T(\overline t_0,\overline y_0;u^{\overline t_0,\overline{y}_0,*}_T)(\overline{t}_0+\delta)) \le (1-\alpha)G_T(\overline{t}_0, \overline{y}_0).
\end{equation}

\end{lemma}
 
 \begin{proof}
 From the definition of $G_T$, in~\eqref{clf}, and Lemma~\ref{lem:bdd_GT}  we have, for~$\delta$ as in Algorithm~\ref{alg:RH} and denoting by simplicity~$\overline\fkY\coloneqq\fkY_T(\overline t_0,\overline y_0;u^{\overline t_0,\overline{y}_0,*}_T)$,
\begin{subequations}\label{e3a}
\begin{align}
&G_T(\overline{t}_0+\delta ,\overline\fkY(\overline{t}_0+\delta)) - G_T(\overline{t}_0, \overline{y}_0)\notag\\
&\hspace{3em}=G_T(\overline{t}_0+\delta ,\overline\fkY(\overline{t}_0+\delta)) -\int^{\overline{t}_0+T}_{\overline{t}_0}\ell^*_T(t;\overline{t}_0, \overline{y}_0)dt-h(T)\norm{\overline{y}_0}{Y}^2\notag \\
&\hspace{3em}\le \left(\Xi(T)- h(T) \right)\norm{\overline{y}_0}{Y}^2+ \left(\Xi(T)- 1 \right) \int_{ \overline{t}_0}^{\overline{t}_0+T} \ell^*_T(t;\overline{t}_0, \overline{y}_0)\rmd t,\\
&\mbox{with}\quad\Xi(T)\coloneqq\tfrac{C_{\rm det}(T)(\gamma(T) + h(T) )}{T-\tau}>0.
\end{align}
\end{subequations}

By using the boundedness of~$\gamma$ and~$C_{\rm det}$ and the properties of~$h$ in~\eqref{clf}, we see that
\[
\lim_{T\to\infty} \Xi(T)=0\qquad\mbox{and}\qquad\lim_{T\to\infty}\tfrac{\Xi(T)}{h(T)} = 0.
\]

Then, for any given~$\alpha\in(0,1)$, we can find a large enough~$T^*>0$, such that
\begin{equation}\label{T(alpha)}
\Xi(T)-1<-\alpha   \qquad\mbox{and}\qquad  \tfrac{\Xi(T)}{h(T)}-1<-\alpha,\qquad\mbox{for all}\quad T\ge T^*.
\end{equation}
This, together with \eqref{e3a}, implies that   
\begin{align}
G_T(\overline{t}_0+\delta ,\overline\fkY(\overline{t}_0+\delta)) - G_T(\overline{t}_0, \overline{y}_0) &\le  - \alpha V_T(\overline{t}_0, \overline{y}_0) -h(T)\alpha \norm{\overline{y}_0}{Y}^2\label{estG-Vic}\\
& \le  -\alpha\int_{\overline{t}_0}^{\overline{t}_0+\delta} \ell^*_T(t;\overline{t}_0, \overline{y}_0)\rmd t,\qquad\mbox{for all}\quad T\ge T^*,\notag
\end{align}
which gives us~\eqref{e20s}.
Finally, by recalling~\eqref{clf}, we observe that~\eqref{estG-Vic} reads 
\begin{equation}
G_T(\overline{t}_0+\delta ,\overline\fkY(\overline{t}_0+\delta))- G_T(\overline{t}_0, \overline{y}_0) \le -\alpha G_T(\overline{t}_0, \overline{y}_0),\notag
\end{equation}
which is equivalent to~\eqref{e20ss}.
 \end{proof}

  For~$T>\tau$, let us denote the level sets within the
  domain~$\clB_r^Y(0)$ of~$G_T(t,\cdot)$ by
   \begin{equation}\label{LevelSets}
  L_{G_{T,t}}(s) \coloneq \{ z \in \clB_r^Y(0) :  G_T(t,z) \leq s  \}, \qquad t\ge0,\quad s\ge0.
   \end{equation}
\begin{lemma}\label{lem:clf}
We have the inequalities
\[ h(T)\norm{\overline{y}_0}{Y}^2 \leq G_T(\overline{t}_0, \overline{y}_0) \leq   (h(T)+\gamma(T))\norm{\overline{y}_0}{Y}^2,   \quad \mbox{for all}\quad (T,\overline{t}_0,\overline{y}_0) \in \bbR_\tau\times\overline\bbR_+ \times \clB_{r}^Y(0).  \]
In particular, $G_T(\overline{t}_0,0)=0$ for all~$\overline{t}_0\ge0$ and~$0\in L_{G_{T,\overline{t}_0}}(s)$ for all~$s\ge0$.
 Furthermore, for every~$0<\rho\le r$, the threshold
$\overline r_\rho(T)\coloneqq\tfrac12h(T)\rho^2>0$
satisfies~$L_{G_{T,\overline t_0}}(\overline r_\rho(T))
\subseteq\clB_\rho^Y(0)$ for all~$\overline t_0\ge0$.
 \end{lemma}
 \begin{proof}
  By using \eqref{clf} and~\eqref{e7} we can see that
 \[
 h(T)\norm{\overline{y}_0}{Y}^2 \le V_T( \overline{t}_0, \overline{y}_0) +h(T)\norm{\overline{y}_0}{Y}^2=G_T( \overline{t}_0, \overline{y}_0)\leq  (\gamma(T)+h(T))\norm{ \overline{y}_0}{Y}^2,
 \]
 which implies~$0\le G_T(\overline t_0,0)\le0$.
 For~$x\in L_{G_{T,\overline t_0}}(\overline r_\rho(T))$,
the lower bound gives
\[
 h(T)\norm{x}{Y}^2
 \le G_T(\overline t_0,x)
 \le\overline r_\rho(T)=\tfrac12h(T)\rho^2.
\]
Since~$h(T)>0$, it follows that
$\norm{x}{Y}\le\rho/\sqrt2<\rho$.
This proves the inclusion in the open ball, uniformly in~$\overline t_0$.
 \end{proof}
\begin{theorem}\label{thm:subopth}
Let Assumption~\ref{asm:rh} hold true. Let~$\alpha \in (0,1)$ and a lower bound $\tau>0$ for sampling time be given. Then there exists $T^* > \tau$ such that for every fixed  prediction horizon $T \geq T^*$,   there exists a radius $ \overline{r} = \overline{r}(T)$ such that for  every  $y_0 \in   L_{G_{T,0}}(\overline{r})$,  the receding-horizon control input~$\mathbf{u}_{\rm rh}$, obtained from  Algorithm \ref{alg:RH},
satisfies
\begin{align}
\label{ed27}
\alpha V_{\infty}(y_0) \leq \alpha J_{\infty}(u_{\rm rh};0,y_0)\leq G_{T}(0,y_0),\\
\label{ed28}
|y_{\rm rh}(t,\cdot )|^2_{Y} \leq C_{\rm rh}e^{-\zeta t}|y_0|^2_{Y} \quad  \text{ for }  t\geq 0,
\end{align}
with $\zeta\coloneqq-T^{-1}\log(1-\alpha)$ and  $C_{\rm rh}\coloneqq\overline C(T)\max \{1, h(T)^{-1}\} \rme^{\zeta  T}(\gamma(T)+h(T))$.
\end{theorem}

\begin{proof}
Since~$\alpha>0$, the definition of~$V_\infty$ implies the first inequality in~\eqref{ed27}. To verify the second one, we choose $T^*=T^*(\tau,\alpha)$ according to Lemma~\ref{lem:estGT}.
For every fixed~$T\ge T^*$, we apply Lemma~\ref{lem:clf} with
$\rho=r_1(T)\le r$ and choose
$\overline r(T)\coloneqq\overline r_{r_1(T)}(T)=\tfrac12h(T)r_1(T)^2$.
Then~$L_{G_{T,t}}(\overline r(T))\subseteq\clB^Y_{r_1(T)}(0)$
for all~$t\ge0$, where~$r_1(T)$ is defined in Lemma~\ref{lem:bdd_GT}.

Next, proceeding by Induction, for the sequence~$(t_k)_{k\in\bbN}$ of concatenation time instants in Algorithm~\ref{alg:RH}  we show that
\begin{align}\stepcounter{equation}
 &y_{\rm rh}(t_{k})\in L_{G_{T,t_{k}}}(\overline{r}),\qquad\mbox{for all}\quad k\in\bbN.\tag{\rm claim.\theequation}\label{k-suboptThm}
\end{align}

\textit{Induction base}: In case~$k=0$, the assumption~$y_0 \in   L_{G_{T,0}}(\overline{r})$ coincides with~\eqref{k-suboptThm}. Therefore, we have that
\begin{equation}\stepcounter{equation}
    \mbox{\eqref{k-suboptThm} holds with $k=0$. }\tag{\rm base.\theequation}\label{kInd_Base}
\end{equation}

\textit{Induction step}: We assume  that for some~$k'\in\bbN$,
\begin{equation}\stepcounter{equation}
    \mbox{\eqref{k-suboptThm} holds with $k=k'$. }\tag{\rm hypo.\theequation}\label{kInd_Hypo}
\end{equation}
Then, since  $y_{\rm rh}(t_{k'}) \in L_{G_{T,t_{k'}}}(\overline{r})$, by  Lemma~\ref{lem:estGT} we have for $\delta_{k'+1} = t_{k'+1}-t_{k'} \ge \tau$ that 
\begin{align}
\label{e88}
&G_T(t_{k'+1},y_{\rm rh}(t_{k'+1})) \leq G_T( t_{k'} ,y_{\rm rh}(t_{k'})) - \alpha \int^{t_{k'+1}}_{t_{k'}} \ell(y_{\rm rh}(t),u_{\rm rh}(t))dt,
\\
\label{e89}
&G_T(t_{k'+1},y_{\rm rh}(t_{k'+1}))  \leq (1-\alpha) G_T(t_{k'},y_{\rm rh}(t_{k'})).
\end{align}
In particular, recalling~\eqref{LevelSets}, we find that
\[
G_T(t_{k'+1},y_{\rm rh}(t_{k'+1}))  \leq (1-\alpha) G_T(t_{k'},y_{\rm rh}(t_{k'}))\le(1-\alpha) \overline r\le\overline r,
\]
which implies that~$y_{\rm rh}(t_{k'+1}) \in L_{G_{T,t_{k'+1}}}(\overline{r})$. Therefore, we have that
\begin{equation}\stepcounter{equation}\label{kInd_Thesis}
    \mbox{\eqref{k-suboptThm} holds with $k=k'+1$. }\tag{\rm thes.\theequation}
\end{equation}
We have just shown that
\begin{equation}\stepcounter{equation}    \mbox{\eqref{kInd_Hypo}\quad$\Longrightarrow$\quad\eqref{kInd_Thesis}}\tag{\rm rec.\theequation}.
\end{equation}
Then, due to~\eqref{kInd_Base}, by induction we can conclude that~\eqref{k-suboptThm} holds true for all~$k\in\bbN$.

Now we observe that, by the arguments above, since~\eqref{k-suboptThm} holds true for all~$k'\in\bbN$, then \eqref{e88}  and \eqref{e89} also hold true for all~$k'\in\bbN$. In particular, by~\eqref{e89} we obtain
\begin{subequations}
\begin{align}
&G_T(t_{k},y_{\rm rh}(t_{k}))  \leq (1-\alpha)^{k-k'} G_T(t_{k'},y_{\rm rh}(t_{k'}))\quad\mbox{for}\quad k\ge k',\label{extGk-geom}
\intertext{and by summing~\eqref{e88} we obtain}
&G_T(t_{k},y_{\rm rh}(t_{k}))- G_T( t_{k'} ,y_{\rm rh}(t_{k'})) \le - \alpha \int^{t_{k}}_{t_{k'}} \ell(y_{\rm rh}(t),u_{\rm rh}(t))\rmd t,\quad\mbox{for}\quad k\ge k'.\label{extGk-arit}
\end{align}
\end{subequations}
Now, taking~$k'=0$ and the limit in~\eqref{extGk-arit} as $k \to \infty$, we find
\[
\lim_{k\to\infty}G_T(t_{k},y_{\rm rh}(t_{k}))- G_T( 0 ,y_{\rm rh}(0)) \le - \alpha \lim_{k\to\infty}\int^{t_{k}}_{0} \ell(y_{\rm rh}(t),u_{\rm rh}(t))\rmd t.
\]
Next, by~\eqref{extGk-geom}, we have that~$G_T(t_{k},y_{\rm rh}(t_{k}))\to0$, which gives us
\begin{equation*}
 \alpha J_{\infty}(u_{\rm rh}, y_0)=\alpha\int^{\infty}_0 \ell(y_{\rm rh}(t),u_{\rm rh}(t))dt \coloneqq\alpha \lim_{k\to\infty}\int^{t_{k}}_{0} \ell(y_{\rm rh}(t),u_{\rm rh}(t)) \rmd t\leq G_T(0,y_0),
\end{equation*}
which is the sought second inequality in~\eqref{ed27}.

It remains to show~\eqref{ed28}. For this purpose we start by observing that for every given $t >0$, we choose the unique integer~$k=k(t) \in \mathbb{N}$ satisfying $t \in [t_k, t_{k+1})$. Using \eqref{e27} we find for each $t \in [t_k, t_{k+1}]\subseteq\clI_{{t}_k,T}$, that 
\begin{align*}  
  \norm{{y}_{\rm rh}(t)}{Y}^2&\le  \overline C(T)\left( \norm{y_{\rm rh}(t_k)}{Y}^2 +\norm{u_T^{t_k,y_{\rm rh}(t_k),*}}{L^2(\clI_{{t}_k,T},\bbR^m)}^2 \right)\\
  &\le  \overline C(T)\left( \norm{y_{\rm rh}(t_k)}{Y}^2 + V_T(t_k,y_{\rm rh}(t_k)) \right)
  \end{align*}
 and, recalling~\eqref{clf} and using~\eqref{e89},
\begin{align*}  
\norm{{y}_{\rm rh}(t)}{Y}^2&\le  \overline C(T)\left( h(T)^{-1}h(T)\norm{y_{\rm rh}(t_k)}{Y}^2 + V_T(t_k,y_{\rm rh}(t_k)) \right)\\
&\le  \overline C(T)\max \{1, h(T)^{-1} \} G_T(t_k,y_{\rm rh}(t_k))\le  C_1(T) (1-\alpha)^kG_T(0,y_{\rm rh}(0)) 
\end{align*}
with~$C_1(T)\coloneqq\overline C(T)\max \{1, h(T)^{-1}\}$.

Now, with~$\zeta\coloneqq -T^{-1}\log(1-\alpha)$, that is, $\rme^{-\zeta T}=(1-\alpha)$,
we obtain
\begin{align*}  
\norm{{y}_{\rm rh}(t)}{Y}^2&\le  C_1(T) \rme^{-\zeta  Tk}G_T(0,y_0) \le  C_1(T) \rme^{\zeta  T}\rme^{-\zeta  t}G_T(0,y_0), 
\end{align*}
because~$t< t_{k+1}\le (k+1)T$. Finally, we use~\eqref{clf} and~\eqref{e7} to find
\begin{align*}  
G_T(0,y_0)=  V_T(0,y_0)+h(T)\norm{{y}_0}{Y}^2 \le (\gamma(T)+h(T))\norm{{y}_0}{Y}^2, 
\end{align*}
which leads us to
\begin{align*}  
\norm{{y}_{\rm rh}(t)}{Y}^2&\le  C_{\rm rh}(T)\rme^{-\zeta  t}\norm{{y}_0}{Y}^2, 
\end{align*}
with~$ C_{\rm rh}(T)\coloneqq \overline C(T)\max \{1, h(T)^{-1}\} \rme^{\zeta  T}(\gamma(T)+h(T))$.
\end{proof}

\begin{remark}
The averaging estimate in the proof of Lemma~\ref{lem:bdd_GT}, together
with~\eqref{e7}, gives
\[
 \norm{y_{\rm rh}(t_{n+1})}{Y}^2
 \le q(T)\norm{y_{\rm rh}(t_n)}{Y}^2,
 \qquad
 q(T)\coloneqq\frac{C_{\rm det}(T)(1+\gamma(T))}{T-\tau}
 \longrightarrow0.
\]
Thus, for sufficiently large $T$, the ball $\clB_r^Y(0)$ is invariant
at the sampling instants, independently of $T$. Since the selected endpoint
remains in this ball, the same Lyapunov-decrease and telescoping arguments
apply on this common neighborhood. For every fixed $y_0\in\clB_r^Y(0)$,
using $V_T(0,y_0)\le V_\infty(y_0)$, we therefore obtain
\[
 V_\infty(y_0)
 \le J_\infty(u_{\rm rh}^{[T]};0,y_0)
 \le\frac{V_\infty(y_0)+h(T)\norm{y_0}{Y}^2}{\alpha(T)}.
\]

 Looking  at the  inequality \eqref{ed27} and~\eqref{T(alpha)}, we can see that as $\alpha\nearrow 1$ and for~$T^*(\alpha)\le T \to \infty$, we will obtain $\lim_{T\to\infty}G_T(0,y_0)\le V_{\infty}(y_0)$. Therefore,  in the limit, 
 \begin{equation}
   V_{\infty}(y_0) = \lim_{T \to \infty}J_{\infty}(u_{\rm rh}^{[T]}, y_0).
 \end{equation}
In this sense, we say that RHC is asymptotically optimal.
\end{remark}


\section{Output-based receding-horizon control for linear systems}
\label{sec:output_feedback}

Here, we assume that the state~$y(t)$ is not fully available. We take~$Y=H=H'$, that is, $H$ is a pivot space. We combine Algorithm~\ref{alg:RH} with a dynamical Luenberger observer to obtain an ORHC (output-based RHC) scheme for the linear system corresponding to
\[
\mathcal{N}(\cdot,y)=A_{\rm rc}(\cdot)\,y,
\]
where~$A_{\rm rc}=A_{\rm rc}(t)\in\clL(V,H)+\clL(H,V')$, with~$V\xhookrightarrow{\rm d,c}H$, is a time-dependent linear reaction--convection-like operator.
Throughout this section, we consider unconstrained controls, that is,
$\bbU=\bbR^m$ (equivalently, $C_u=\infty$). Thus, on every finite time
interval $I$, $L^2_{\rm ad}(I,\bbR^m)=L^2(I,\bbR^m)$.
Note that, now, we cannot proceed with Step~\ref{alg:RH-update-state} of Algorithm~\ref{alg:RH}, because the full state is not available. Instead, we employ sensor measurements \(w(t)\in\bbR^{S}\) to compute an estimate \(\widehat y(t_n)\approx y(t_n)\) at the sampling instants. The control input is then implemented as
\[
u(t)=u(\widehat y(t_n);t),\qquad t\in(t_n,t_{n+1}),
\]
which yields an ORHC framework. The ORHC state--input pair \((y_{\rm orh},u_{\rm orh})\) satisfies
\begin{align}\label{systRHCtildeu_short}
\dot y_{\rm orh}+Ay_{\rm orh}+A_{\rm rc}y_{\rm orh}=B u_{\rm orh},\qquad
w_{\rm orh}=\clZ y_{\rm orh},
\end{align}
with \(y_{\rm orh}(0)=y_0\), and
\[
u_{\rm orh}\rest{(t_n,t_{n+1})}
=u_T^*(t;t_n,\widehat y(t_n))\rest{(t_n,t_{n+1})},\qquad n\in\bbN.
\]
The measured output is given by
\begin{equation}\label{wZy_short}
w_{\rm orh}(t)=\clZ y_{\rm orh}(t)\coloneqq(\Psi_1 y_{\rm orh}(t),\dots,\Psi_S y_{\rm orh}(t))\in\bbR^S.
\end{equation}

The estimate \(\widehat y\) is generated by the full-dimensional observer
\begin{align}\label{eq:obsTV_noq}
\dot{\widehat y}+A\widehat y+A_{\rm rc}\widehat y
=Bu_{\rm orh}+\fkI_S(\clZ\widehat y-w_{\rm orh}),\qquad
\widehat y(0)=\widehat y_0,
\end{align}
where \(\fkI_S(t):\bbR^S\to H\) denotes the output-injection operator. We choose/construct \(\fkI_S\) such that, for all \((\widehat y_0,y_0)\in H\times H\),
\begin{equation}\label{goal-obs_short}
\norm{\widehat y(t)-y_{\rm orh}(t)}{H} \le C_1 e^{-\mu_1 t}\norm{\widehat y_0-y_0}{H},\qquad t\ge0,
\end{equation}
with a constant \(C_1\) that is independent of \(u_{\rm orh}\) (simply because the estimation error dynamics is independent of the input, due to the linearity of the dynamics). The ORHC procedure is summarized in Algorithm~\ref{alg:ORHC}.

\begin{algorithm}[htbp]
\caption{Output-based receding horizon control, ORHC($\tau,T;\clZ,\fkI_S,\widehat y_0$)}\label{alg:ORHC}
\begin{algorithmic}[1]
\REQUIRE{lower bound $\tau$ for sampling time; prediction horizon $T\ge \tau$; output operator $\clZ$; injection operator $\fkI_S$; initial estimate $\widehat y_0\in H$.}
\ENSURE{a control input $u_{\rm orh}\in L^2_{\rm ad}(\bbR_+,\bbR^m)$; observer/state estimate $\widehat y(\cdot)$; sampling times $(t_k)_{k\in\bbN}$.}

\STATE Set $(\overline t_0,\overline{y}_0):=(0,\widehat y_0)$, $t_0:=0$, and initialize the observer with $\widehat{y}(0)=\widehat y_0$.
\STATE Set $T_\infty:=\infty$; $k:=0$.
\WHILE{$\overline t_0<T_\infty$}
  \STATE \label{alg:ORHC_solve}``Solve'' the  FTH problem~\ref{tag:opfin} with $\mathcal{N}(\cdot, y) = A_{\rm rc}(\cdot)y$ to obtain an optimal input
  $u_T^{\overline t_0,\overline{y}_0,*}\in L^2(\clI_{\overline t_0,T},U)$ and the associated predicted trajectory
  $y_T^{\overline{t}_0,\overline{y}_0,*}(\cdot)=\fkY_T(\overline t_0,\overline{y}_0;u_T^{\overline{t}_0,\overline{y}_0,*})(\cdot)$.
  \STATE \label{alg:ORHC_delta} Set
  $\delta := \max\argmin_{s\in[\tau,T]} \norm{y_T^{\overline t_0,\overline{ y}_0,*}(\overline t_0+s)}{H}$.
  \STATE Apply the control on $[\overline t_0,\overline t_0+\delta)$:
  \[
  u_{\rm orh}(t):=u_T^{\overline t_0,\overline{ y}_0,*}(t),\qquad t\in[\overline t_0,\overline t_0+\delta).
  \]
   
  \STATE Let the plant \eqref{systRHCtildeu_short} evolve on $[\overline t_0,\overline t_0+\delta)$ with input $u_{\rm orh}$ and record the output
  $w_{\rm orh}(t):=\clZ y_{\rm orh}(t)$.
  
  \STATE Propagate the observer on $[\overline t_0,\overline t_0+\delta)$ using the recorded output:
  \[
  \dot{\widehat y}+A\widehat y+A_{\rm rc}\widehat y
  =B u_{\rm orh}+\fkI_S(\clZ\widehat y-w_{\rm orh}),
  \qquad \widehat y(\overline t_0)=\overline{y}_0.
  \]
  
  \STATE Update:
  $(\overline{t}_0,\overline{y}_0)\leftarrow (\overline t_0+\delta,\widehat y(\overline t_0+\delta))$.
  \STATE $k\gets k+1$;\quad $t_k\gets \overline t_0$.
\ENDWHILE
\end{algorithmic}
\end{algorithm}

\begin{lemma}[Lipschitz bounds for LQR on a fixed horizon]\label{lem:lipschitz-lqr-T}
Assume the finite-horizon problem on $\clI_{\overline{t}_0,T}$ is linear--quadratic with fixed horizon $T>0$. 
For each initial time $\overline{t}_0\ge0$, the finite-horizon value function
is represented by a bounded, self-adjoint, nonnegative operator
$P_T(\overline{t}_0)$, namely
\[
V_T(\overline{t}_0,y)=\langle P_T(\overline{t}_0)y,y\rangle_H,\qquad y\in H.
\]
Moreover, by~\ref{asm:rh-bddVT}, we have
$\norm{P_T(\overline{t}_0)}{\mathcal L(H)}\le\gamma(T)$
uniformly in $\overline{t}_0\ge0$. The
optimal control is given by a linear state feedback along the optimal trajectory. In particular, there exist constants $L_p = L_p(T),L_u=L_u(T)\ge0$ such that for all $\overline{t}_0\ge0$ and all $y,z\in H$,
\begin{equation}
\label{eq:Lip_cont_lqr}
\norm{P_T({\overline{t}_0})}{\mathcal L(H)} \leq L_p, \quad \text{and} \quad  \norm{u_T^*(\cdot;\overline{t}_0,y)-u_T^*(\cdot;\overline{t}_0,z)}{L^2(\clI_{\overline{t}_0,T},\bbR^m)}
\le L_u(T)\,\norm{y-z}{H},
\end{equation}
In particular, $G_T$ has a quadratic form. That is, $G_T(\overline{t}_0,y)=\langle (P_T(\overline{t}_0)+h(T)\mathcal{I})y,y\rangle_H$, where $\mathcal{I}$ stands for the identity operator. Moreover,  for all $y,z\in H$ and all $\varepsilon>0$, we have
\begin{equation}\label{eq:quad-pert-lemma}
\begin{split}
G_T(\overline{t}_0,y+z)\le (1+\varepsilon)\,G_T(\overline{t}_0,y)+\Bigl(1+\frac{1}{\varepsilon}\Bigr)\,G_T(\overline{t}_0,z).
\end{split}
\end{equation}
\end{lemma}

\begin{proof}
The first parts of the  proof follow by standard arguments. To show \eqref{eq:quad-pert-lemma}, for simplicity of the presentation we set
$G \coloneq G_T(\overline{t}_0,\cdot)$ and $P \coloneq P_T(\overline{t}_0)+h(T)\mathcal{I}$. Then, using the fact that  $G(y+z)=\langle P(y+z),y+z\rangle_H$, we expand
\[
G(y+z)=G(y)+2\langle Py,z\rangle_H+G(z).
\]
Since $P\ge0$, we may write $\langle Py,z\rangle_H=\langle P^{1/2}y,P^{1/2}z\rangle_H$. Hence, by Cauchy--Schwarz and
Young's inequality, for any $\varepsilon>0$,
\[
\begin{split}
&2|\langle Py,z\rangle_H|
=2|\langle P^{1/2}y,P^{1/2}z\rangle_H|
\le 2 \norm{P^{1/2}y}{H}\,\norm{P^{1/2}z}{H}
 \\ &\le \varepsilon\norm{P^{1/2}y}{H}^2+\frac{1}{\varepsilon}\norm{P^{1/2}z}{H}^2
=\varepsilon G(y)+\frac{1}{\varepsilon}G(z).
\end{split}
\]
Substituting into the expansion of $G(y+z)$ yields
\[
G(y+z)\le (1+\varepsilon)G(y)+\Bigl(1+\frac{1}{\varepsilon}\Bigr)G(z),
\]
which proves \eqref{eq:quad-pert-lemma}.
\end{proof}

\begin{theorem}\label{Thm:orhc-noq}
Let Assumption~\ref{asm:rh} hold and fix $\alpha\in(0,1)$ and $\tau>0$.
Let $T\ge T^*$ and $\overline r=\overline r(T)$ be as in Theorem~\ref{thm:subopth}.
Consider ORHC (Algorithm~\ref{alg:ORHC}) for the time-varying plant \eqref{systRHCtildeu_short} with observer~\eqref{eq:obsTV_noq}. Denote by $(y_{\rm orh},u_{\rm orh})$ the ORHC closed-loop pair and set $e(t)\coloneqq \widehat y(t)-y_{\rm orh}(t)$,
$e_0\coloneqq \widehat y_0-y_0$. Assume that there exist $C_1\ge1$ and $\mu_1>0$ such that
\begin{equation}
\label{eq:exponent_output_inj}
\norm{e(t)}{H}\le C_1 e^{-\mu_1 t}\norm{e_0}{H},\qquad t\ge 0.
\end{equation}
Then, for every $\varepsilon\in\bigl(0,\tfrac{\alpha}{1-\alpha}\bigr)$,
 for the given system and chosen observer, there exist constants
$C_{\rm orh}\ge1$ and $\mu_2>0$, depending on $T,\tau,\alpha,\varepsilon$
but independent of the initial conditions $y_0$ and $\widehat y_0$,
such that the ORHC trajectory satisfies
\begin{equation}\label{eq:exp-decay-y}
\norm{y_{\rm orh}(t)}{H}^2\le C_{orh} \,e^{-\mu_2 t }\Big(\norm{\widehat{y}_0}{H}^2+\norm{e_0}{H}^2\Big),\qquad t \geq 0,
\end{equation}
\end{theorem}

\begin{proof}
For a given $n\in\bbN$, we set
$\delta_n:=t_{n+1}-t_n\in[\tau,T]$ and denote the
optimal prediction actually used by Algorithm~\ref{alg:ORHC} by
\[
p_n(t):=y_T^*(t;t_n,\widehat y(t_n)),
\qquad p_n(t_n)=\widehat y(t_n).
\]
On $[t_n,t_{n+1}]$, the plant and $p_n$ receive the same input
$u_{\rm orh}=u_T^*(\cdot;t_n,\widehat y(t_n))$. Consequently, their difference
$d_n:=y_{\rm orh}-p_n$ satisfies
\[
\dot d_n+Ad_n+A_{\rm rc}(t)d_n=0,
\qquad d_n(t_n)=-e(t_n).
\]
By homogeneity, the estimates in Assumption~\ref{asm:rh} and~\eqref{e20ss}
hold globally in the present unconstrained linear--quadratic setting.
Applying~\eqref{e27} to the homogeneous equation with zero input, extended
to $[t_n,t_n+T]$, gives
\begin{equation}
\label{eq:z_e_estimate}
\norm{d_n}{C([t_n,t_{n+1}];H)}
\le D(T)\norm{e(t_n)}{H},
\qquad D(T):=\sqrt{\overline C(T)}.
\end{equation}
Since $\delta_n$ is selected from $p_n$ itself, the nominal
decrease~\eqref{e20ss}, time-shifted to $t_n$, applies to this prediction:
\[
G_T(t_{n+1},p_n(t_{n+1}))
\le (1-\alpha)G_T(t_n,\widehat y(t_n)).
\]
Set $C_G(T):=\gamma(T)+h(T)$ and $G_n:=G_T(t_n,y_{\rm orh}(t_n))$.
Then Lemma~\ref{lem:clf} gives $G_T(t,\xi)\le C_G(T)\norm{\xi}{H}^2$.
For the prescribed $\varepsilon\in(0,\alpha/(1-\alpha))$, we choose
\[
\eta:=\sqrt{1+\varepsilon}-1>0,
\qquad (1+\eta)^2=1+\varepsilon.
\]
Since $G_T(t,\xi)=\langle(P_T(t)+h(T)\mathcal I)\xi,\xi\rangle$ is
quadratic, we can apply~\eqref{eq:quad-pert-lemma} twice, first at the
endpoint and then to $\widehat y(t_n)=y_{\rm orh}(t_n)+e(t_n)$. This yields
\begin{align*}
G_{n+1}
&\le (1+\eta)G_T(t_{n+1},p_n(t_{n+1}))
 +(1+\eta^{-1})C_G(T)D(T)^2\norm{e(t_n)}{H}^2\\
&\le (1+\eta)(1-\alpha)G_T(t_n,\widehat y(t_n))
 +(1+\eta^{-1})C_G(T)D(T)^2\norm{e(t_n)}{H}^2\\
&\le (1+\eta)^2(1-\alpha)G_n
 +(1+\eta^{-1})C_G(T)
 \bigl[D(T)^2+(1+\eta)(1-\alpha)\bigr]\norm{e(t_n)}{H}^2.
\end{align*}
Thus, with $\rho:=(1-\alpha)(1+\varepsilon)<1$, we obtain
\begin{equation}\label{eq:GT-rec2-short}
G_{n+1}\le \rho\,G_n+\beta\norm{e(t_n)}{H}^2,
\qquad \beta:=(1+\eta^{-1})C_G(T)
\bigl[D(T)^2+(1+\eta)(1-\alpha)\bigr].
\end{equation}
Invoking the observer estimate $\norm{e(t_n)}{H}\le C_1 e^{-\mu_1 t_n}\norm{e_0}{H}$ gives
\[
G_{n+1}\le \rho\,G_n + d\,e^{-2\mu_1 t_n}\norm{e_0}{H}^2,
\qquad d:=\beta C_1^2 .
\]
Iterating this recursion yields, for all $n\ge1$, 
\begin{equation}\label{eq:GT-iter-short}
G_n \le \rho^n G_0
+ d\sum_{j=0}^{n-1}\rho^{n-1-j}e^{-2\mu_1 t_j}\|e_0\|_H^2 .
\end{equation}
Since $t_j\ge j\tau$, we have $e^{-2\mu_1 t_j}\le e^{-2\mu_1 j\tau}$, hence
\begin{equation}\label{eq:GT-iter2-short}
G_n \le \rho^n G_0
+ d\sum_{j=0}^{n-1}\rho^{n-1-j}e^{-2\mu_1 j\tau}\|e_0\|_H^2.
\end{equation}
Setting $r:=e^{-2\mu_1\tau}/\rho$ and
\[
S_{n-1}:=\sum_{j=0}^{n-1}\rho^{n-1-j}e^{-2\mu_1 j\tau}
=\rho^{n-1}\sum_{j=0}^{n-1}r^j .
\]
we have the following cases: 
\[
S_{n-1}\le 
\begin{cases}
\frac{1}{1-r}\,\rho^{n-1}, & 0<r<1,\\[0.6ex]
\frac{r}{r-1}\,e^{-2\mu_1 (n-1)\tau}, & r>1,\\[0.6ex]
C_\theta\,\rho^{\theta (n-1)}, & r=1,\ \theta\in(0,1),
\end{cases}
\]
where in the case $r=1$ we used $S_{n-1}=n\rho^{n-1}$ and the estimate
$n\rho^{n-1}\le C_\theta \rho^{\theta(n-1)}$ for $n\ge1$ and the choice $C_\theta:=\frac{\rho^{\theta-1}}{e(1-\theta)|\log\rho|}$.
Consequently, we have
\begin{equation}\label{eq:Gn-unique-mu}
G_n \le C_{\rm rec}\,e^{-\mu n\tau}\bigl(G_0+\|e_0\|_H^2\bigr),
\end{equation}
where, for any fixed $\theta\in(0,1)$ in the case $r=1$, we choose
$C_{\rm rec}:=\max\{1,d A_r e^{\mu\tau}\}$, with
\begin{equation}\label{eq:mu-cases}
\mu:=
\begin{cases}
\min\Bigl\{\dfrac{|\log\rho|}{\tau},\,2\mu_1\Bigr\}, & r\neq 1,\\[1.2ex]
\min\Bigl\{\dfrac{\theta|\log\rho|}{\tau},\,2\mu_1\Bigr\}, & r=1,
\end{cases}
\qquad
A_r:=\begin{cases}
(1-r)^{-1},&0<r<1,\\
r/(r-1),&r>1,\\
C_\theta,&r=1.
\end{cases}
\end{equation}
Here $C_\theta$ bounds the sum in the case $r=1$, whereas
$C_{\rm rec}$ is the constant in~\eqref{eq:Gn-unique-mu}; when $r=1$,
$C_{\rm rec}$ also depends on the chosen $\theta$ through $A_r=C_\theta$.
Since $C_{\rm rec}\ge1$, estimate~\eqref{eq:Gn-unique-mu} also holds for $n=0$.
Using $G_0\le C_G(T)\norm{y_0}{H}^2$,
$y_0=\widehat y_0-e_0$, and the lower bound in Lemma~\ref{lem:clf}, we obtain
\[
\norm{y_{\rm orh}(t_n)}{H}^2 \le \widehat C(T)\,e^{-\mu n \tau }\bigl(\norm{\widehat y_0}{H}^2+\norm{e_0}{H}^2\bigr),\qquad n\in\bbN,
\]
where $\widehat C(T):=C_{\rm rec}(2C_G(T)+1)/h(T)$; dependence on the other fixed
parameters is suppressed in this notation.

Now let $t\ge0$ and choose $n\in\bbN$ with $t\in[t_n,t_{n+1}]$.
The predicted optimal input satisfies~\eqref{contC} with initial state
$\widehat y(t_n)$, so~\eqref{e27} implies
\[
\norm{p_n}{C([t_n,t_n+T];H)}^2
\le \overline C(T)(1+\gamma(T))\norm{\widehat y(t_n)}{H}^2.
\]
Using $y_{\rm orh}=p_n+d_n$, estimate~\eqref{eq:z_e_estimate}, and
$\widehat y(t_n)=y_{\rm orh}(t_n)+e(t_n)$, we therefore obtain
\begin{align}
\norm{ y_{\rm orh}(t)}{H}^2
&\le 2\overline C(T)(1+\gamma(T))\norm{\widehat y(t_n)}{H}^2
 +2D(T)^2\norm{e(t_n)}{H}^2\notag\\
&\le C_*(T)\left(\norm{y_{\rm orh}(t_n)}{H}^2
 +\norm{e(t_n)}{H}^2\right),\label{eq:thm_orhs_exponent0}
\end{align}
where $C_*(T):=4\overline C(T)(1+\gamma(T))+2D(T)^2$.
Since $t_n\ge n\tau$ and $\mu\le2\mu_1$, the observer estimate gives
\[
\norm{e(t_n)}{H}^2
\le C_1^2 e^{-2\mu_1t_n}\norm{e_0}{H}^2
\le C_1^2 e^{-\mu n\tau}\norm{e_0}{H}^2.
\]
Combining this with the sampling-time bound, for every
$n\in\bbN$, including $n=0$, we~have
\[
\norm{y_{\rm orh}(t)}{H}^2
\le C_2(T)e^{-\mu n\tau}
\bigl(\norm{\widehat y_0}{H}^2+\norm{e_0}{H}^2\bigr),
\qquad C_2(T):=C_*(T)(\widehat C(T)+C_1^2).
\]
In particular, no identification of $d_n$ with the observer error is used.
Finally, $t\le t_{n+1}\le(n+1)T$ implies
$n\tau\ge(\tau/T)t-\tau$, and hence
$e^{-\mu n\tau}\le e^{\mu\tau}e^{-(\mu\tau/T)t}$.
Thus~\eqref{eq:exp-decay-y} follows with
$C_{\rm orh}:=\max\{1,C_2(T)e^{\mu\tau}\}$ and
$\mu_2:=\mu\tau/T$.
This completes the proof.
\end{proof}

\begin{remark}
Note that we are particularly interested in the case~$Q$ has finite-dimensional range. However, we would like to mention that if we choose $Q$, in~\eqref{Obj}, as the identity operator from $Y$ onto $Y$, then the conclusions of Theorem~\ref{Thm:orhc-noq} remain valid without the need of squeezing norm step~\ref{alg:ORHC_delta} in Algorithm~\ref{alg:ORHC}, and with an a priori fixed $\tau>0$ and time instances $t_i=i\tau$. In this case, we can work with $G_T(t,z)\coloneqq V_T(t,z)$ (cf.~\eqref{clf-GT=VT+h}). By \cite[Prop.~2.8]{AzmiKartVolk26}, this function serves as a Lyapunov functional and decreases exponentially to zero. Moreover, by invoking \cite[Thm.~2.9]{AzmiKartVolk26}, the exponential stability and suboptimality estimates \eqref{ed27}--\eqref{ed28} hold for $V_T$ in place of $G_T$.
 Using also the endpoint estimate in~\cite[Lem.~2.5, eq.~(5)]{AzmiKartVolk26}, estimates analogous to those in the proof of Theorem~\ref{Thm:orhc-noq} follow with $V_T$ in place of $G_T$.
We therefore emphasize that the necessity (within our proof) of the squeezing norm step~\ref{alg:ORHC_delta} in Algorithm~\ref{alg:ORHC} arises from allowing $Q$ in $\ell$ to have a nontrivial kernel.
\end{remark}

\section{Applicability of the abstract RHC assumptions to parabolic equations}
\label{sec:parabolic_verification}

In this section we verify Assumption~\ref{asm:rh} for three representative
parabolic models under the hypotheses stated below. The first example is the Schl\"ogl equation with
finite-dimensional internal control. The second is the Schl\"ogl equation
with dynamic finite-dimensional boundary control, written as an extended
system. The third is a nonautonomous linear parabolic equation with the same
dynamic boundary-control mechanism; this is the model used in the
output-feedback setting of Section~\ref{sec:output_feedback}.

For each model the verification follows the same three steps. Requirement
\ref{asm:rh-existopt} concerns existence of an optimal control for every
finite-horizon problem. Requirement~\ref{asm:rh-bddVT} follows from a
finite-dimensional stabilizing feedback and the resulting uniform upper bound
for the value function. Requirement~\ref{asm:rh-estimYT} combines
well-posedness and a priori estimates with the detectability-type estimate
\eqref{e8}, which recovers the full state from the semidefinite state cost and
the control effort. We make these three steps explicit below and retain the
detailed estimates needed for the less standard~\ref{asm:rh-estimYT} verification.

The proof mechanisms are common to the three examples. For~\ref{asm:rh-existopt} we use the
direct method of the calculus of variations: a minimizing sequence is bounded
by the coercive input term, the corresponding states satisfy uniform energy
estimates, and weak compactness together with sequential weak lower
semicontinuity yields a minimizer. For the Schl\"ogl nonlinearity, the required
continuity estimates are available in~\cite[Sect.~III]{AzmiKunRod23} and
\cite[Sect.~3.2]{KunRod23-dcds}. For the extended boundary-controlled system,
the complete optimal-control existence argument is given in
Subsection~\ref{ssec:Schl-bdry}, while the supporting well-posedness result is
proved in Appendix~\ref{app:schl-wellposedness}.

For~\ref{asm:rh-bddVT}, we first write the estimate for a time-independent
linear feedback, as in the internal constructions
in~\cite{KunRod23-dcds,Rod21-aut}: $u^m(t)=\fkK_my^m(t)$, where
$\fkK_m\in\clL(Y,\bbR^m)$ and
$\kappa_m=\norm{\fkK_m}{\clL(Y,\bbR^m)}<\infty$.
Suppose this feedback is admissible and that the associated state satisfies
$\norm{y^m(\overline t_0+s)}{Y}
\le C_{\rm st}\rme^{-\omega s}\norm{\overline y_0}{Y}$ for $s\ge0$.
Here $C_{\rm st}\ge1$ and $\omega>0$ may depend on $m$, but are uniform
over $\overline t_0\ge0$ and $\overline y_0\in\clB_r^Y(0)$.
With $q=\norm{Q}{\clL(Y,\bbR^s)}$, insertion into~\eqref{Obj} gives
\begin{align}
 V_T(\overline t_0,\overline y_0)
 &\le J_T(u^m;\overline t_0,\overline y_0)\le C_{\rm st}^2(q^2+\kappa_m^2)
       \norm{\overline y_0}{Y}^2\int_0^T\rme^{-2\omega s}\,\rmd s
   =\gamma(T)\norm{\overline y_0}{Y}^2,\label{eq:rh2-feedback-cost}\\
 \gamma(T)&\coloneqq
       \frac{C_{\rm st}^2(q^2+\kappa_m^2)}{2\omega}
                         (1-\rme^{-2\omega T}).\label{eq:rh2-gamma}
\end{align}
Thus $\gamma$ is continuous, nondecreasing, and bounded, as required in~\ref{asm:rh-bddVT},
for each fixed $m$.
The same estimate holds for a time-dependent feedback if $\kappa_m$ is a
uniform bound for its operator norm; see~\cite{PhanRod18-mcss} for related
boundary-feedback constructions.

In the boundary-control examples, $\fkK_m$ acts on the
extended state $y=(z-\clE\clB v,v)$ and produces the
auxiliary input $u$, whereas the boundary input control forcing is $\clB v$.
Here $q$ is understood as the norm of the bounded state-cost
map $y\mapsto(Q(y_1+\clE\clB y_2),y_2)$.

For~\ref{asm:rh-estimYT}, estimate~\eqref{e27} follows from the standard
energy estimate for the corresponding state equation. The essential point is
therefore~\eqref{e8}: the finite-dimensional operator~$Q$ must be rich enough
to control the full state norm. This detectability-type estimate is verified
explicitly in each example.

\subsection{Example I: Schl\"ogl equation with internal controls and internal state-cost operators}
\label{ssec:Sch-internal}
We consider the equation
 \begin{align}\label{Schlogl}
    &\tfrac{\p}{\p t} y=\nu\Delta y-y-\clP(t,y)+Bu,\quad  (\bfn\cdot\nabla y)\rest{\p\Omega}=0,\quad y(0)=y_0
\end{align}
which takes the form~\eqref{cs}, with free dynamics  defined by the operators
\begin{equation}
\label{operSchlogl}
A=A_{0}\coloneqq-\nu\Delta+\Id,\qquad\clN(t,y)=\clP(t,y)\coloneqq y(y-2)(y+2)-y.
\end{equation}
Here we consider~$H\coloneqq L^2(\Omega)$ as the pivot space~$H=H'$. The operator~$A_{0}\colon V\to V'$ stands for the shifted Neumann Laplacian, with the Sobolev space~$V=W^{1,2}(\Omega)\subset H$. In other words, we consider here homogeneous Neumann boundary conditions~$(\bfn\cdot\nabla y)\rest{\p\Omega}=0$ and
\[
\langle A_0 y,h\rangle_{V',V}\coloneqq\nu(\nabla y,\nabla h)_{L^2(\Omega)^d}+(y,h)_{L^2(\Omega)},\qquad (y,h)\in V\times V.
\]

The well-posedness is known for strong solutions. Thus, we consider initial states~$y_0\in V$.

The bounded spatial domain~$\Omega\subset\bbR^d$, $d\in\{1,2,3\}$, is either a smooth domain or a convex polygonal/polyhedral domain.

For simplicity, we denote the pivot space as~$H\coloneqq L^2(\Omega)$, and we recall that we also have
$A_0\colon \rmD(A_0)\to H$ with domain
\begin{equation}\label{DNeuLap}
 \rmD(A_0)=\{h\in W^{2,2}(\Omega)\mid(\bfn\cdot\nabla h)\rest{\p\Omega}=0\}. 
\end{equation}

\subsubsection{Verification of~\ref{asm:rh-existopt}}
For every finite horizon and every initial state in~$Y=V$, existence and
uniqueness of the strong state follow from the standard Schl\"ogl estimates.
The input term in the cost bounds every minimizing sequence in
$L^2(I;\bbR^m)$. The compactness and continuity properties established
in~\cite[Sect.~III]{AzmiKunRod23} and
\cite[Sect.~3.2]{KunRod23-dcds}, followed by weak lower semicontinuity of
the cost, give an optimal control. Hence~\ref{asm:rh-existopt} holds for the internal-control
model.

\subsubsection{Verification of~\ref{asm:rh-bddVT} and admissible finite-dimensional configurations}
For sufficiently many localized internal actuators, finite-dimensional
stabilizing feedbacks are available by the arguments of
\cite{KunRod23-dcds}. The feedback produces a uniformly bounded
infinite-horizon cost and therefore the~\ref{asm:rh-bddVT} estimate. We record below the
state-cost operators used in the subsequent~\ref{asm:rh-estimYT} verification.

\begin{figure}[htbp]%
\centering%
\includegraphics[width=1\textwidth]{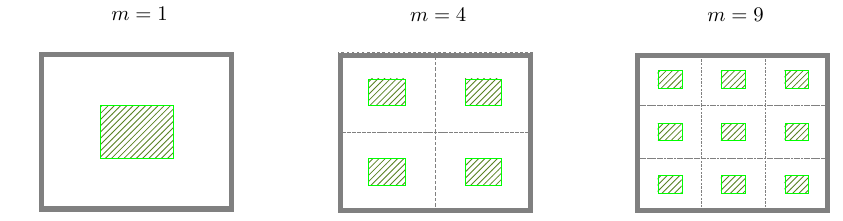}
\caption{Supports of actuators on a rectangular domain~$\Omega$.} \label{fig:suppActGrid}
\end{figure}

\paragraph{Eigenfunction-based state-cost operator.}
We consider the case~$B\colon\bbR^m\to H$ and~$Q\colon H\to \bbR^s$ are given as follows:
\begin{equation}
\label{operBQ}
Bu=\sum_{i=1}^m u_i\indf_{\omega_i},\qquad Qy=\left((y,e_1)_H,(y,e_2)_H,\dots,(y,e_s)_H\right)%
\end{equation}
where~$\indf_{\omega_i}$ are the indicator functions of subdomains~${\omega_i}$ as in Fig.~\ref{fig:suppActGrid}
and the~$e_j$ are ``the'' first orthonormalized eigenfunctions of~$A$ (ordered by increasing eigenvalues repeated according to their multiplicity).
The satisfiability of~\eqref{e8} in Assumption~\ref{asm:rh} follows by the results in~\cite{KunRod23-dcds}. See, in particular~\cite[Sect.~3.1]{KunRod23-dcds}.
\begin{remark}
In fact, in~\cite{KunRod23-dcds} the integrand of the cost functional was taken as
$
\widehat\ell(y,u)=\norm{P_{\clE_s^\rmf}y}{H}^2+\norm{u}{\bbR^m}^2
$
instead of
$
\ell(y,u)=\norm{Qy}{\bbR^s}^2+\norm{u}{\bbR^m}^2
$
as in~\eqref{Obj}. However, we can see that these are equivalent for an orthonormal system of eigenfunctions.
\end{remark}

\begin{figure}[htbp]%
\centering%
\includegraphics[width=1\textwidth]{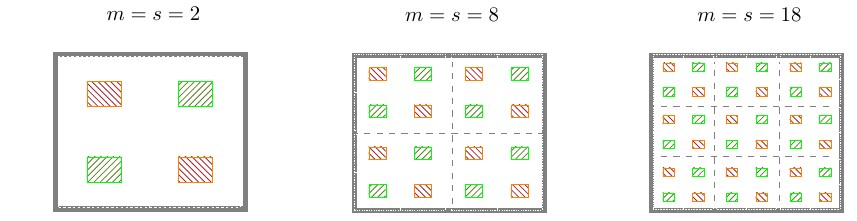}
\caption{Supports of actuators (``$\slash$''-pattern) and of ``sensors'' (``$\backslash$''-pattern).} \label{fig:suppActSensGrid}
\end{figure}

\paragraph{Localized state-cost operator.}
We consider the case~$B\colon\bbR^m\to H$ and~$Q\colon H\to \bbR^s$ are given as follows:
\begin{equation}
\label{operBQ-new}
Bu=\sum_{i=1}^{m} u_i\indf_{\omega^1_i},\qquad Qy=\left((y,\overline\indf_{\omega^2_1})_H,(y,\overline\indf_{\omega^2_2})_H,\dots,(y,\overline\indf_{\omega^2_s})_H\right)%
\end{equation}
where the~$\omega^1_i$ and the~${\omega^2_i}$ are the domains as in Fig.~\ref{fig:suppActSensGrid}

A motivation for considering normalized sensor indicator functions~$\overline\indf_{\omega^2_j}=\norm{\indf_{\omega^2_j}}{H}^{-1}\indf_{\omega^2_j}$ instead of eigenfunctions~$e_j$ is that the former involve local ``measurements'' only and do not require the computations of the eigenfunctions of~$A$. Here, we could indeed see~$Qy$ as the output of local average-like measurements performed by sensors~$\indf_{\omega^2_j}$. In this case, the optimal control goal would be motivated by the stabilization of such output. However, in this manuscript,~$Qy$ does not need to correspond to an output of concrete sensor measurements.

Global stabilizability, for a sufficiently large number~$m$ of actuators,
with control operators as in~\eqref{operBQ-new} follows by arguments analogous
to those in~\cite{KunRod23-dcds} and completes~\ref{asm:rh-bddVT} for this configuration.

\subsubsection{Verification of~\ref{asm:rh-estimYT}}
We verify the detectability estimate~\eqref{e8} with $Y=V$ for the
localized state-cost operator in~\eqref{operBQ-new}, adapting the energy
argument of~\cite{KunRod23-dcds}. We retain the regularized indicator functions
\begin{equation}\label{reg-act}
 \phi_j(x)=\indf_{\omega_j^2}(x)\prod_{k=1}^d
 \sin\!\left(\frac{\pi(x_k-a_{jk})}{b_{jk}-a_{jk}}\right),
 \qquad \omega_j^2\coloneqq\bigtimes_{k=1}^d(a_{jk},b_{jk}),
\end{equation}
and set
\[
 \clW_s\coloneqq\linspan\{\indf_{\omega_j^2}:1\le j\le s\},
 \qquad
 \widetilde\clW_s\coloneqq\linspan\{\phi_j:1\le j\le s\}\subset V.
\]
The sensor supports are pairwise disjoint, and
$(\phi_i,\indf_{\omega_j^2})_H=0$ for $i\ne j$, whereas
$(\phi_j,\indf_{\omega_j^2})_H>0$. Hence we conclude that
$H=\widetilde\clW_s\oplus\clW_s^\perp$, and we write
\[
 \Pi_s\coloneqq P_{\widetilde\clW_s}^{\clW_s^\perp},\qquad
 P_s\coloneqq P_{\clW_s}^{\clW_s^\perp},\qquad
 \theta\coloneqq\Pi_s y,\qquad
 \Theta\coloneqq(\Id-\Pi_s)y
       =P_{\clW_s^\perp}^{\widetilde\clW_s}y.
\]
Thus $y=\theta+\Theta$, with $\theta\in\widetilde\clW_s$ and
$\Theta\in V\cap\clW_s^\perp$ whenever $y\in V$.
Since $\Pi_s=\Pi_sP_s$ and its range is a finite-dimensional subspace
of $V$, the constants
$c_{s,H}\coloneqq\norm{\Pi_s}{\clL(H)}$ and
$c_{s,V}\coloneqq\norm{\Pi_s}{\clL(H,V)}$
are finite, and we have
\[
 \norm{\theta}{H}\le c_{s,H}\norm{P_s y}{H},\qquad
 \norm{\theta}{V}\le c_{s,V}\norm{P_s y}{H}.
\]
For the spatially distributed refining sensor family in
Fig.~\ref{fig:suppActSensGrid}, the Poincar\'e-type estimate
in~\cite[Sect.~6]{Rod21-aut} gives
\[
 \beta_s\coloneqq
 \inf_{\Theta\in(V\cap\clW_s^\perp)\setminus\{0\}}
       \frac{\norm{\Theta}{V}^2}{\norm{\Theta}{H}^2}
 \longrightarrow+\infty.
\]
For this construction $c_{s,H}$ is also uniformly bounded, but no
uniform bound on $c_{s,V}$ is needed below.
The normalized sensor indicators in~\eqref{operBQ-new} form an
orthonormal basis of $\clW_s$, so that
$\norm{P_s y}{H}^2=\norm{Qy}{\bbR^s}^2$. Consequently, we can write
\begin{equation}\label{eq:localized-sensor-estimate}
\begin{aligned}
 \norm{y}{H}^2
 &\le 2\norm{\theta}{H}^2+2\norm{\Theta}{H}^2
 \le 2c_{s,H}^2\norm{P_s y}{H}^2
                   +2\beta_s^{-1}\norm{y-\theta}{V}^2\\
 &\le 4\beta_s^{-1}\norm{y}{V}^2
       +D_s\norm{Qy}{\bbR^s}^2,
 \qquad D_s\coloneqq2c_{s,H}^2+4\beta_s^{-1}c_{s,V}^2.
\end{aligned}
\end{equation}

Let $I_T=(\overline t_0,\overline t_0+T)$, with $0<T<\infty$, and
write $y_0=y(\overline t_0)$. Testing~\eqref{Schlogl} with $y$ and
using $\clP(y)=y^3-5y$ yields
\[
 \frac12\frac{\rmd}{\rmd t}\norm{y}{H}^2
   +\norm{y}{V}^2+\norm{y}{L^4(\Omega)}^4
 =5\norm{y}{H}^2+(Bu,y)_H.
\]
Young's inequality and time integration therefore give
\[
 \norm{y}{L^2(I_T,V)}^2
 \le C_1\left(\norm{y_0}{H}^2
          +\norm{y}{L^2(I_T,H)}^2
          +\norm{Bu}{L^2(I_T,H)}^2\right),
\]
where $C_1>0$ is independent of $\overline t_0$, $T$, and $s$.
Substituting the preceding sensor estimate, we obtain
\[
 \left(1-\frac{4C_1}{\beta_s}\right)
                    \norm{y}{L^2(I_T,V)}^2
 \le C_1\Bigl(\norm{y_0}{H}^2
          +D_s\norm{Qy}{L^2(I_T,\bbR^s)}^2
          +\norm{Bu}{L^2(I_T,H)}^2\Bigr).
\]
Choose $s$ along the above refining family so that $\beta_s\ge8C_1$
and then keep this sensor configuration fixed. The coefficient on the
left is at least $1/2$. Since $B\in\clL(\bbR^m,H)$, it follows that
\begin{align}\label{eq:localized-rh3-bound}
 \norm{y}{L^2(I_T,V)}^2
 &\le C_2\left(\norm{y_0}{H}^2
          +\norm{Qy}{L^2(I_T,\bbR^s)}^2
          +\norm{u}{L^2(I_T,\bbR^m)}^2\right)\notag\\
 &=C_2\left(\norm{y_0}{H}^2
          +J_T(u;\overline t_0,y_0)\right)
 \le C_2\left(\norm{y_0}{V}^2
          +J_T(u;\overline t_0,y_0)\right),
\end{align}
where
$C_2\coloneqq2C_1\max\{1,D_s,\norm{B}{\clL(\bbR^m,H)}^2\}$
is independent of $\overline t_0$, $T$, $y_0$, and $u$, but may depend
on the fixed actuator and sensor configurations. This proves~\eqref{e8}
with $Y=V$ and $C_{\rm det}(T)=C_2$. The same estimate holds for
$T=\infty$ by passage to the limit from finite intervals whenever the
right-hand side is finite.

\begin{remark}\label{rem:uniform-sensor-bound}
It is worth mentioning that, for the configurations in
Fig.~\ref{fig:suppActSensGrid}, the constant $D_s$
in~\eqref{eq:localized-sensor-estimate} is uniformly bounded.
Indeed, $c_{s,H}$ is uniformly bounded. Writing $s=2N^2$ for these
configurations, where $N$ is the refinement index, the scaling calculation
in~\cite{Rod21-aut} and the Poincar\'e argument in~\cite{Rod21-sicon} give
$c_{s,V}^2\le c_{s,H}^2(1+a_1N^2)$ and $\beta_s\ge1+a_2N^2$,
with $a_1,a_2>0$ independent of $N$. The factor $c_{s,H}^2$ in the
first bound comes from applying the norm comparison on
$\widetilde\clW_s$ to $\Pi_s y$. Thus $\beta_s^{-1}c_{s,V}^2$, and
hence $D_s$, are uniformly bounded. Since the actuator norms are also
uniformly bounded, $C_2$ in~\eqref{eq:localized-rh3-bound} can be chosen
independently of $s$, for all sufficiently fine configurations.
\end{remark}

\begin{remark}
    Above, in Fig.~\ref{fig:suppActGrid}, for the sake of a simpler illustration, we consider a rectangular spatial domain~$\Omega$. The strategy can be applied to more general domains; see~\cite[Rem.~2.8]{AzmiKunRod23}.
\end{remark}

\paragraph{Coincident actuator and state-cost supports.}
With the actuators distributed as in Fig.~\ref{fig:suppActGrid}, for large~$m$, we can also take
 \begin{equation}
\label{operBQ-same}
Bu=\sum_{i=1}^m u_i\indf_{\omega_i},\qquad Qy=\left((y,\overline\indf_{\omega_1})_H,(y,\overline\indf_{\omega_2})_H,\dots,(y,\overline\indf_{\omega_m})_H\right).%
\end{equation}

\paragraph{Actuation in a prescribed subdomain.}
Instead of actuators given by indicator functions~$\indf_{\omega_j}$ as in  Fig.~\ref{fig:suppActGrid} or~$\indf_{\omega_j^1}$ as in Fig.~\ref{fig:suppActSensGrid} we can also take indicator functions all supported in an a-priori given subdomain. The only difference is that, in general, we may have local stabilizability only. Such local stabilizability holds for a sufficiently large number~$m$ of actuators; see the arguments in~\cite{PhanRod18-mcss}. 

\subsection{Example II: Schl\"ogl equation with dynamic boundary control and internal state-cost operators}\label{ssec:Schl-bdry}
We first rewrite the boundary-controlled equation as an extended system with a
bounded control operator. We then verify~\ref{asm:rh-existopt}--\ref{asm:rh-estimYT} for that formulation, with~\ref{asm:rh-bddVT} conditional on the stabilizability hypothesis stated below.

The Schl\"ogl equation, under dynamic Neumann boundary controls, with actuators~$\psi_j\in W^{\frac12,2}(\p\Omega)$, $1\le j\le m$, assumed to be linearly independent, reads as follows
\begin{subequations}\label{Schlogl-ext0}
 \begin{align}\label{Schlogl-1}
    &\tfrac{\p}{\p t} z=\nu\Delta z-z-\clP(\cdot,z),\quad  (\bfn\cdot\nabla z)\rest{\p\Omega}=\clB v\coloneqq\sum_{j=1}^m v_j\psi_{j},\quad&& z(0)=z_0,\\
    &\tfrac{\rmd}{\rmd t} v=-v+u,\quad&&v(0)=v_0.\label{Schlogl-2}
\end{align}
\end{subequations}
which we can write in the form of system~\eqref{cs},
\begin{subequations}\label{Schlogl-ext}
\begin{align}
&\dot y+Ay+\clN(\cdot,y)=Bu,\qquad y(0)=y_0;
\intertext{with state~$y=(y_1,y_2)\coloneqq  (z-\clE \clB v,v)\in V\times\bbR^m$, $A_0$ and~$\clP$ as in~\eqref{operSchlogl}, $V=\rmD(A_0^\frac12)=W^{1,2}(\Omega)$, and with}
&Ay=\begin{bmatrix}A_0&0\\
0&\Id\end{bmatrix}y,\qquad  Bu=\begin{bmatrix}-\clE \clB\\\Id\end{bmatrix}u,\\
&\clN(\cdot,y)=\begin{bmatrix}\clP(\cdot,y_1+\clE \clB y_2)+\overline A_0\clE \clB y_2-\clE\clB y_2\\0\end{bmatrix}.\label{Schlogl-ext-N}
\end{align}
\end{subequations}

Here, the pivot space is taken as~$H\coloneqq L^2(\Omega)\times\bbR^m$, and the initial conditions will be taken as~$y_0=(y_{1,0},y_{2,0})\in Y \coloneqq V \times\bbR^m$.
We can see that the operator~$A$ has domain~$\rmD(A)=\rmD(A_0)\times\bbR^m$, with~$\rmD(A_0)$ as in~\eqref{DNeuLap}. 
 Above, the~$\psi_j\in W^{\frac12,2}(\p\Omega)$, $1\le j\le m$, are the Neumann boundary actuators and~$\clE $ is an extension operator lifting the actuators as follows
 \begin{equation}\label{liftE}
    \Psi_j\coloneqq \clE \psi_j\in W^{2,2}(\Omega)\qquad\mbox{with}\qquad(\bfn\cdot\nabla\Psi_j)\rest{\p\Omega}=\psi_j.
 \end{equation}
 The operator~$\overline A_0$ is defined in the range of~$\clE $ by
\begin{align*}
\langle \overline A_0\clE \psi_j,h\rangle_{V',V}&\coloneqq\langle (-\nu\Delta +\Id)(\clE \psi_j),h\rangle_{V',V}\\
&\coloneqq\langle A_0(\clE \psi_j),h\rangle_{V',V}-\nu(\psi_j,h)_{L^2(\p\Omega)},\quad\mbox{for all}\quad h\in V.
\end{align*}

 Note that, if~$(\bfn\cdot\nabla z)\rest{\p\Omega}=\clB v$, then~$(\bfn\cdot\nabla y_1)\rest{\p\Omega}=(\bfn\cdot\nabla (z-\clE \clB v))\rest{\p\Omega}=0$.
Note also that for~$(y_1,y_2)\in\rmD(A)$ we have that
\begin{align*}
&\langle\overline A_0(y_1+\clE \clB y_2),h\rangle_{V',V}=\langle A_0y_1+\overline A_0\clE \clB y_2,h\rangle_{V',V}\\
&\hspace{4em}=\langle A_0(y_1+\clE \clB y_2),h\rangle_{V',V}-\nu(\clB y_2,h)_{L^2(\p\Omega)}=\langle(-\nu\Delta+\Id)(y_1+\clE \clB y_2),h\rangle_{V',V}.
\end{align*}
That is,~$\overline A_0$ coincides with the formal Laplacian~$-\nu\Delta+\Id$ in the space
$\rmD(A_0)\oplus \clE \clB\bbR^m$, which is isomorphic to~$\rmD(A_0)\times\bbR^m=\rmD(A)$.

 We can take, for example, the operator~$\clE =\clE ^0$ defined as follows:
 \[
(-\nu\Delta +\Id)(\clE ^0\psi_j)=0,\quad (\bfn\cdot\nabla(\clE ^0\psi_j))\rest{\p\Omega}=\psi_j,\qquad 1\le j\le m,
 \]
leading to
\[
\langle \overline A_0\clE ^0\psi_j,h\rangle_{V',V}
=\langle (-\nu\Delta +\Id)(\clE ^0\psi_j),h\rangle_{L^2(\Omega)}
=0,\quad\mbox{for all}\quad h\in V.
\]
Thus this particular lifting eliminates the interior term
$\overline A_0\clE ^0\psi_j$ from~\eqref{Schlogl-ext-N}; the boundary action
is already encoded by the decomposition $z=y_1+\clE ^0\clB v$.

\subsubsection{Verification of~\ref{asm:rh-existopt}}
\label{ssec:exist-optcontrol-Schlogl}

The well-posedness result needed below is established in
Appendix~\ref{app:schl-wellposedness}. In particular, for every bounded
interval~$I=(\overline t_0,\overline t_0+T)$, every initial state
$y_0\in Y=V\times\bbR^m$, and every input $u\in L^2(I;\bbR^m)$, the extended
system~\eqref{Schlogl-ext} has a unique strong solution
\[
 y^u\in W(I;\rmD(A),H)\xhookrightarrow{}C(\overline I;Y).
\]
We now use this result to verify existence of an optimal control, which is the
substantive assertion in~\ref{asm:rh-existopt}.

\begin{proposition}[Existence of a finite-horizon optimal control]
\label{prop:schl-bdry-optimal-control}
Let $I=(\overline t_0,\overline t_0+T)$ be a nonempty bounded interval and let
$y_0\in Y$. Suppose that the pointwise admissible set $\bbU\subseteq\bbR^m$ is
nonempty, closed, and convex; the unconstrained case $\bbU=\bbR^m$ is included.
Then the finite-horizon problem~\ref{tag:opfin}, with the state-cost operator
and extended state described in~\eqref{J-Qz.v.u}, admits an optimal control.
\end{proposition}

\begin{proof}
Fix $y_0=(y_{10},y_{20})\in V\times\bbR^m$ and
$I=(\overline t_0,\overline t_0+T)$, with $\overline t_0\ge0$ and $T>0$.
Let $(u^k)_{k\in\bbN}\subset L^2_{\rm ad}(I;\bbR^m)$ be a minimizing
sequence, that is,
\[
 \clJ(u^k)\longrightarrow
       \inf_{u\in L^2_{\rm ad}(I;\bbR^m)}\clJ(u),
 \qquad \clJ(u)\coloneqq J(y^u,u).
\]
Here $J$ is the functional in~\eqref{J-Qz.v.u}; explicitly,
\begin{align*}
 \clJ(u)&=\norm{Q(y_1^u+\clE\clB y_2^u)}{L^2(I;\bbR^s)}^2+\norm{y_2^u}{L^2(I;\bbR^m)}^2
                   +\norm{u}{L^2(I;\bbR^m)}^2,
\end{align*}
where $y^u=(y_1^u,y_2^u)$ is the solution of~\eqref{Schlogl-ext} on~$I$
with initial state~$y_0$.
Write $y^k=y^{u^k}=(y_1^k,y_2^k)$ and
$z^k=y_1^k+\clE\clB y_2^k$. Since the cost contains
$\norm{u}{L^2(I;\bbR^m)}^2$, the sequence $(u^k)$ is bounded in
$L^2(I;\bbR^m)$. After passage to a subsequence,
\[
 u^k\xrightharpoonup[L^2(I;\bbR^m)]{}u^*.
\]
The admissible set is closed and convex, hence weakly closed, and therefore
$u^*\in L^2_{\rm ad}(I;\bbR^m)$.

The uniform estimate~\eqref{Galerkin-unifbound} in
Theorem~\ref{thm:schl-wellposedness} bounds the sequence~$(y^k)$ in
$W(I;\rmD(A),H)$. Consequently, after taking a further
subsequence,
\begin{align*}
 y^k&\xrightharpoonup[L^2(I;\rmD(A))]{}y^*,
 &\dot y^k&\xrightharpoonup[L^2(I;H)]{}\dot y^*,
 &y^k&\xrightarrow[L^2(I;Y)]{}y^*.
\end{align*}
The last convergence follows from the Aubin--Lions lemma because
$\rmD(A)\xhookrightarrow{\rm c}Y\xhookrightarrow{}H$. In particular,
$z^k\to z^*\coloneqq y_1^*+\clE\clB y_2^*$ in $L^2(I;V)$. The polynomial
identity used in the proof of Theorem~\ref{thm:schl-wellposedness}, together
with the uniform $L^\infty(I;V)$ bound, gives
\[
 \clP(\cdot,z^k)\xrightharpoonup[L^2(I;L^2(\Omega))]{}
 \clP(\cdot,z^*).
\]
We may therefore pass to the limit in the weak formulation of
\eqref{Schlogl-ext}. The continuity of the trace map on
$W(I;\rmD(A),H)$ preserves the initial condition, and uniqueness yields
$y^*=y^{u^*}$.

Finally, the state-cost operator is continuous on~$Y$. Hence the two state
terms in~\eqref{J-Qz.v.u} converge along the strongly convergent state
sequence, while the squared control norm is weakly lower semicontinuous.
Therefore,
\[
 J(y^{u^*},u^*)\leq\liminf_{k\to\infty}J(y^{u^k},u^k)
 =\inf_{u\in L^2_{\rm ad}(I;\bbR^m)}J(y^u,u).
\]
Thus $u^*$ is optimal:
$\clJ(u^*)=\min_{u\in L^2_{\rm ad}(I;\bbR^m)}\clJ(u)$.
\end{proof}

Proposition~\ref{prop:schl-bdry-optimal-control} proves~\ref{asm:rh-existopt} for
$Y=V\times\bbR^m$.

\subsubsection{Verification of~\ref{asm:rh-bddVT}}
For finite-dimensional Dirichlet boundary controls, related stabilization
results are established in~\cite{PhanRod18-mcss}; see also
\cite[Thm.~1.6 and Sect.~3.2]{BadraTakahashi14} for spectral stabilization
and dynamic-extension arguments. For the dynamic Neumann boundary-control
system considered here, we assume that the extended system~\eqref{Schlogl-ext} admits an
admissible bounded linear feedback that locally exponentially stabilizes
the origin in $Y=V\times\bbR^m$, with stability constants independent of
the initial time. The existence of such a feedback is an additional
hypothesis and is not established in this paper. Under this hypothesis,
inserting the feedback into the cost functional, as
in~\eqref{eq:rh2-feedback-cost}, yields~\ref{asm:rh-bddVT} locally.
An analogous uniform stabilizability hypothesis is retained for the
nonautonomous linear Neumann example in Subsection~\ref{ssec:linear-bdry-obs}.

\subsubsection{Verification of~\ref{asm:rh-estimYT}}
The well-posedness part of~\ref{asm:rh-estimYT} and estimate~\eqref{e27} follow from the strong
solution bounds established in Appendix~\ref{app:schl-wellposedness}. It remains
to verify the detectability estimate~\eqref{e8}. The following energy argument
does so for sufficiently many localized internal state-cost functions.

We distinguish the variational operator~$A_0:V\to V'$ from the
differential operator~$\overline A_0$ acting on states with nonhomogeneous
Neumann data. We have
\[
 \langle A_0z,z\rangle_{V',V}=\norm{z}{V}^2,
 \qquad
 \langle\overline A_0z,z\rangle_{V',V}
 =\norm{z}{V}^2-\nu(\clB v,z)_{L^2(\p\Omega)}.
\]
Thus the boundary flux is counted only once. Testing~\eqref{Schlogl-ext0}
with~$(z,v)$, and setting
$E(t)\coloneqq\norm{z(t)}{L^2(\Omega)}^2+\norm{v(t)}{\bbR^m}^2$, gives
\begin{equation}\label{energy-bdry1}
\begin{split}
 E'(t)
 &=-2\langle\overline A_0z,z\rangle_{V',V}
   -2(\clP(\cdot,z),z)_{L^2(\Omega)}
   -2\norm{v}{\bbR^m}^2+2(u,v)_{\bbR^m}\\
 &=-2\norm{z}{V}^2+2\nu(\clB v,z)_{L^2(\p\Omega)}
   -2(\clP(\cdot,z),z)_{L^2(\Omega)}
   -2\norm{v}{\bbR^m}^2+2(u,v)_{\bbR^m}.
\end{split}
\end{equation}
For the polynomial~$\clP(\cdot,z)=z(z-2)(z+2)-z=z^3-5z$, we obtain
\begin{align}
 -2(\clP(\cdot,z),z)_{L^2(\Omega)}
 &=-2\norm{z}{L^4(\Omega)}^4+10\norm{z}{L^2(\Omega)}^2
 \le C\norm{z}{L^2(\Omega)}^2,\qquad C=10.
 \label{energy-P-Schlogl}
\end{align}
The trace theorem and boundedness of the finite-dimensional boundary
operator give a constant~$C_{\partial}\ge0$ such that
\[
 \nu\big|(\clB v,z)_{L^2(\p\Omega)}\big|
 \le C_{\partial}\norm{v}{\bbR^m}\norm{z}{V}.
\]
Consequently, we can write
\begin{equation}\label{energy-bdry2}
 E'(t)\le-2\norm{z}{V}^2
   +2C_{\partial}\norm{v}{\bbR^m}\norm{z}{V}
   +C\norm{z}{L^2(\Omega)}^2
   -2\norm{v}{\bbR^m}^2+2(u,v)_{\bbR^m}.
\end{equation}
Using Young's inequality in the forms
\[
 2C_{\partial}\norm{v}{\bbR^m}\norm{z}{V}
 \le\norm{z}{V}^2+C_{\partial}^2\norm{v}{\bbR^m}^2,
 \qquad
 2(u,v)_{\bbR^m}\le\tfrac12\norm{u}{\bbR^m}^2+2\norm{v}{\bbR^m}^2,
\]
and retaining the damping term~$-2\norm{v}{\bbR^m}^2$ until it cancels,
we find
\begin{align}
 E'(t)+\norm{z}{V}^2
 &\le C_{\partial}^2\norm{v}{\bbR^m}^2
   +C\norm{z}{L^2(\Omega)}^2+\tfrac12\norm{u}{\bbR^m}^2.
 \label{est-schl-bdry-energy}
\end{align}
For the localized sensor family in~\eqref{operBQ-same}, with supports as
in Fig.~\ref{fig:suppActGrid}, the coercivity result of~\cite{AzmiKunRod23}
gives, for every~$\mu>0$, sufficiently large~$s$ and~$\lambda>0$ such that
\[
 \norm{z}{V}^2+2\lambda\norm{P_{W_s}z}{L^2(\Omega)}^2
 \ge\mu\norm{z}{L^2(\Omega)}^2,\qquad z\in V.
\]
Here~$P_{W_s}$ is the orthogonal projection in~$L^2(\Omega)$.
Since the sensors are orthonormal,
$\norm{P_{W_s}z}{L^2(\Omega)}^2=\norm{Qz}{\bbR^s}^2$.
Choosing~$\mu=2C=20$ and absorbing half of the~$V$-norm yields
\[
 E'(t)+\tfrac12\norm{z}{V}^2
 \le C_{\partial}^2\norm{v}{\bbR^m}^2
    +\lambda\norm{Qz}{\bbR^s}^2+\tfrac12\norm{u}{\bbR^m}^2.
\]
For any~$\overline t_0\ge0$ and~$T>0$, set
$I=(\overline t_0,\overline t_0+T)$. Time integration gives
\[
 E(\overline t_0+T)+\tfrac12\norm{z}{L^2(I;V)}^2
 \le E(\overline t_0)+C_{\partial}^2\norm{v}{L^2(I;\bbR^m)}^2
       +\lambda\norm{Qz}{L^2(I;\bbR^s)}^2
       +\tfrac12\norm{u}{L^2(I;\bbR^m)}^2.
\]
Dropping the nonnegative final energy and then adding the integral of
$\norm{v}{\bbR^m}^2$, we obtain
\begin{align*}
 \norm{(z,v)}{L^2(I;Y)}^2
 &\le2E(\overline t_0)+(2C_{\partial}^2+1)\norm{v}{L^2(I;\bbR^m)}^2
       +2\lambda\norm{Qz}{L^2(I;\bbR^s)}^2
       +\norm{u}{L^2(I;\bbR^m)}^2\\
 &\le C_2\big(E(\overline t_0)+J(y,u)\big),
\end{align*}
where~$C_2\coloneqq\max\{2,2C_{\partial}^2+1,2\lambda\}$ and~$J$ is
given in~\eqref{J-Qz.v.u} below.

To pass to the extended state, put~$L\coloneqq\clE\clB$ and
$C_L\coloneqq\max\{2,1+2\norm{L}{\mathcal L(\bbR^m,V)}^2\}$.
The squared triangle inequality, with its factor~$2$, gives
\[
 \norm{y}{Y}^2
 =\norm{z-Lv}{V}^2+\norm{v}{\bbR^m}^2
 \le2\norm{z}{V}^2+2\norm{Lv}{V}^2+\norm{v}{\bbR^m}^2
 \le C_L\norm{(z,v)}{Y}^2.
\]
Similarly, using~$\norm{h}{L^2(\Omega)}\le\norm{h}{V}$, the inverse
change of variables at the initial time satisfies
\begin{align*}
 E(\overline t_0)
 &=\norm{y_1(\overline t_0)+Lv(\overline t_0)}{L^2(\Omega)}^2
    +\norm{v(\overline t_0)}{\bbR^m}^2\\
 &\le2\norm{y_1(\overline t_0)}{L^2(\Omega)}^2
    +2\norm{Lv(\overline t_0)}{L^2(\Omega)}^2
    +\norm{v(\overline t_0)}{\bbR^m}^2
 \le C_L\norm{y(\overline t_0)}{H}^2.
\end{align*}
Combining these estimates and taking~$C_5=C_L^2C_2$ and~$C_6=C_5$
proves~\eqref{e8}:
\begin{align}
 \norm{y}{L^2(I;Y)}^2
 &\le C_5\left(\norm{y(\overline t_0)}{H}^2
       +\norm{Qz}{L^2(I;\bbR^s)}^2
       +\norm{v}{L^2(I;\bbR^m)}^2
       +\norm{u}{L^2(I;\bbR^m)}^2\right)
       \label{obsEstim-Schl-bdry-pivoty0}\\
 &\le C_6\left(\norm{y(\overline t_0)}{Y}^2+J(y,u)\right),
       \label{obsEstim-Schl-bdry}
\end{align}
with the unchanged cost
\begin{equation}\label{J-Qz.v.u}
J(y,u)\coloneqq\norm{Q(y_1+\clE\clB y_2)}{L^2(I;\bbR^s)}^2+\norm{y_2}{L^2(I;\bbR^m)}^2+\norm{u}{L^2(I;\bbR^m)}^2.
\end{equation}
The constant~$C_6$ is independent of~$(\overline t_0,T)$; it may depend
on the fixed diffusion coefficient, domain, actuators, sensors, and lifting
operator.

\begin{remark}
By considering the extended system~\eqref{Schlogl-ext} above we look for the input~$v$ in a dynamic form; \cite{PhanRod18-mcss}.
See also~\cite[Sect.~4.2]{BadraTakahashi11} in an autonomous setting. Another approach we find in the literature is to look for the input~$v$ in a static form, that is, to omit~\eqref{Schlogl-2} and take~$v$ as the control variable in~\eqref{Schlogl-1}, ~\cite{BarLasTri06,Raymond06}. Note that we cannot apply our theory directly to this case, because~\eqref{Schlogl-1} alone is not in the form of~\eqref{cs}.
\end{remark}

\subsection{Example III: Linear parabolic equation with dynamic boundary control and output feedback}\label{ssec:linear-bdry-obs}
We consider a linear system of the form~\eqref{systRHCtildeu_short} and an
observer of the form~\eqref{eq:obsTV_noq}, as in
Section~\ref{sec:output_feedback}, now combined with dynamic boundary control.

For example, we can consider the equation
\begin{subequations}\label{sys-z-lin-bdry}
 \begin{align}
    &\tfrac{\p}{\p t} z=\nu\Delta z-z-\clN(t,z),\qquad  (\bfn\cdot\nabla z)\rest{\p\Omega}=\clB v,\qquad z(0)=z_0, \\
    &\dot v=-v+u,\qquad v(0)=v_0
\end{align}    
\end{subequations}
evolving in a bounded spatial domain~$\Omega\subset\bbR^d$, $d\in\bbN_+$, which takes the form~\eqref{cs} once we consider the analogous extended system~\eqref{Schlogl-ext}. The reaction-convection operator is taken as
\begin{equation}
\label{Arc-lin-oper-bdry}
\clN(t,z)=A_{\rm rc}(t)z\coloneqq az+b\cdot\nabla z,
\end{equation}
with suitable coefficients~$a=a(t,x)\in\bbR$ and~$b=b(t,x)\in\bbR^d$, assumed to be essentially bounded, $a\in L^\infty(\bbR_+\times\Omega)$ and~$b\in (L^\infty(\bbR_+\times\Omega))^d$.
Here we use the pivot and state space
$H=Y\coloneqq L^2(\Omega)\times\bbR^m$, with~$H=H'$. As in Example~II,
the extended operator contains the shifted Neumann Laplacian
$A_0=-\nu\Delta+\Id$ in its first component and the identity in its second
component.

In this case an observer for~\eqref{sys-z-lin-bdry} can be designed as (cf.~\eqref{eq:obsTV_noq})
\begin{subequations}\label{sys-zhat-lin-bdry}
 \begin{align}
    &\tfrac{\p}{\p t} \widehat z=\nu\Delta \widehat z-\widehat z-\clN(t,\widehat z)+\fkI_S(Q\widehat z-Q z),\qquad  (\bfn\cdot\nabla \widehat z)\rest{\p\Omega}=\clB v,\qquad \widehat z(0)=\widehat z_0, \\
    &\dot v=-v+u,\qquad v(0)=v_0.
\end{align}    
\end{subequations}

\subsubsection{Verification of~\ref{asm:rh-existopt}}
For fixed initial data, the linear state equation defines a continuous affine
control-to-state map. The coercive quadratic input term bounds every minimizing
sequence, and the direct method for linear-quadratic problems yields an optimal
control; see, for example,~\cite[Thm.~1.43]{HinPinUlb08}. Thus~\ref{asm:rh-existopt} holds
globally in the pivot space~$Y=H$.

\subsubsection{Verification of~\ref{asm:rh-bddVT}}
Assume that the extended linear Neumann system admits a uniformly bounded
linear stabilizing feedback, with exponential stability constants independent
of the initial time. Inserting this feedback into the cost functional, as
in~\eqref{eq:rh2-feedback-cost}, yields the required uniform bound for the
value function. Thus,~\ref{asm:rh-bddVT} holds under this stabilizability
hypothesis. For related boundary-feedback constructions, see, for
example,~\cite{Rod21-amo,PhanRod18-mcss}.

\subsubsection{Verification of~\ref{asm:rh-estimYT}}
Linear parabolic well-posedness and the standard energy estimate imply
\eqref{e27}. To establish~\eqref{e8}, with state-cost functions as
in~\eqref{operBQ-same} and supports as in Fig.~\ref{fig:suppActGrid}, we
proceed as in Example~II. Energy estimates for~\eqref{sys-z-lin-bdry} give the
analogue of~\eqref{energy-bdry1},
\begin{align*}
    \tfrac{\rmd}{\rmd t} \norm{(z,v)}{L^2(\Omega)\times\bbR^m}^2
    &=-2\norm{z}{V}^2
      +2\nu(\clB v,z)_{L^2(\p\Omega)}
      -2(A_{\rm rc}(\cdot,z),z)_{L^2(\Omega)}\\
    &\quad-2\norm{v}{\bbR^m}^2+2(u,v)_{\bbR^m}.
\end{align*}
For the linear term we will have, using~\cite[Lem.~3.5]{Rod21-aut}, that for any~$\gamma>0$,
\begin{align*}
 -2(A_{\rm rc}(\cdot)z,z)_{L^2(\Omega)} &\le 2\gamma \norm{z}{V}^2+2\overline C\gamma^{-1}\norm{z}{L^2(\Omega)}^2,
\end{align*}
for some constant~$\overline C\ge0$. Taking~$\gamma=\frac1{4}$, we obtain the slight variation of~\eqref{energy-bdry2},
\begin{align*}
    &\tfrac{\rmd}{\rmd t} \norm{(z,v)}{L^2(\Omega)\times\bbR^m}^2\le-\norm{z}{V}^2+2C_1\norm{v}{\bbR^m}\norm{z}{W^{1,2}(\Omega)}+C\norm{z}{L^2(\Omega)}^2+2(u,v)_{\bbR^m}.
\end{align*}
with~$C=2\overline C\gamma^{-1}=8\overline C$. Now we can proceed with a small variation of the arguments in Section~\ref{ssec:Schl-bdry} to arrive at the analogue of~\eqref{obsEstim-Schl-bdry}, which gives us
     \begin{align*}
  \norm{y}{L^2(I;H)}^2 \le  \norm{y}{L^2(I;V\times\bbR^m)}^2
     &\le\widetilde C_{6}\left(\norm{y(\overline t_0)}{H}^2+J(y,u)\right),
\end{align*}
again with~$J(y,u)\coloneqq\norm{Q(y_1+\clE \clB y_2)}{L^2(I;\bbR^s)}^2+\norm{y_2}{L^2(I;\bbR^m)}^2+\norm{u}{L^2(I;\bbR^m)}^2$. This gives us the sought~\eqref{e8}, with a constant~$\widetilde C_6$ independent of the pair~$(t_0,T)$.

\begin{remark}
    Recall that $z=y_1+\clE \clB y_2$ is the physical state. Thus the cost
    functional~$J(y,u)$ penalizes finite-dimensional measurements of the
    physical state. The detectability estimate~\eqref{e8} then transfers the
    decay of these measurements and of the input terms to the full state.
\end{remark}

For the output-feedback result, one additionally chooses the injection
operator~$\fkI_S$ so that the observer error satisfies
\eqref{eq:exponent_output_inj}. This requirement concerns the observer and is
separate from~\ref{asm:rh-existopt}--RH3, which have just been verified for the controlled plant.

\subsection{Scope and choice of the state space}
The preceding examples concern Schl\"ogl and linear parabolic equations, but
analogous arguments apply to wider classes of semilinear parabolic evolution
equations.

These examples also show that the state space~$Y$ depends on the model. When
strong solutions are available but weak solutions for arbitrary
$y_0\in H=H'$ are not, the initial condition is taken in
$Y\subset L^2(\Omega)=H$ for internal controls, or in
$Y\subset L^2(\Omega)\times\bbR^m=H$ for dynamic boundary controls; see
Sections~\ref{ssec:Sch-internal} and~\ref{ssec:Schl-bdry}. If weak
well-posedness is known, it can be convenient to take $Y=H$; see
\cite[Rem.~1.3]{KunRodWal21} and Sections~\ref{sec:output_feedback}
and~\ref{ssec:linear-bdry-obs}. If the actuators are distributions, such as
Dirac deltas, it may instead be convenient to take $Y=V'\supset H$ and study
stabilizability in the weaker $V'$-norm; see~\cite{KunRodWal24}.

\section{Numerical experiments}
\label{sec:numerics}
We consider two cases: the Schl\"ogl equation with boundary state-feedback controls, and a linear model with dynamic output-feedback controls. For numerical experiments in the case of a linear model with internal controls we refer the reader to~\cite{Azmi25output}.

\begin{figure}[!htbp]
  \centering
  \includegraphics[width=\textwidth]{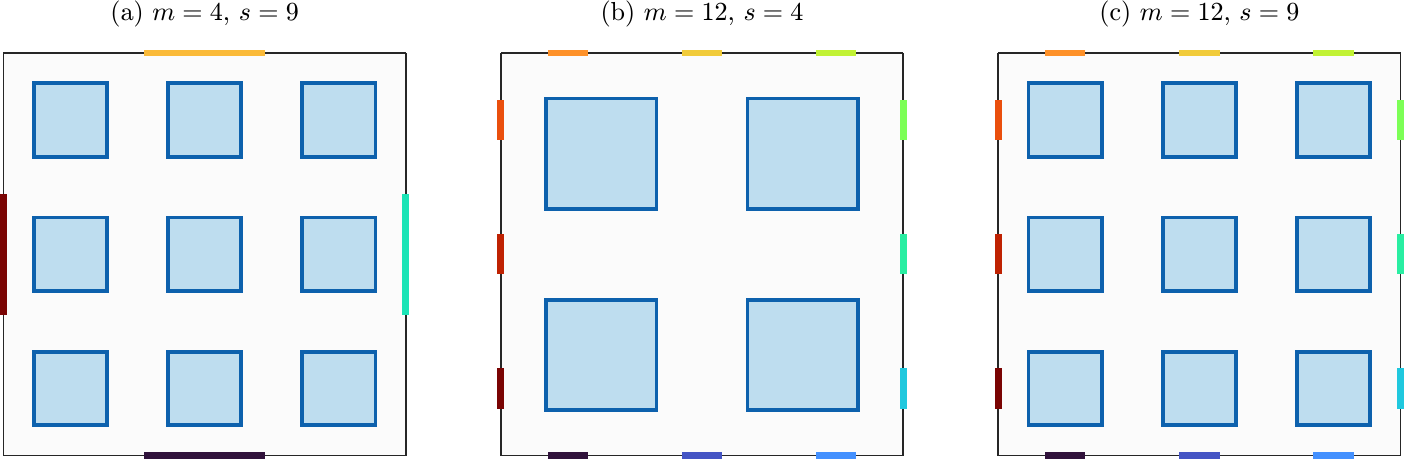}
  \caption{Supports of~$m$ boundary actuators and of~$s$ internal sensors.}
  \label{fig:placements}
\end{figure}

The spatial domain is taken as the unit square $\Omega=(0,1)\times(0,1)$. The diffusion parameter is taken as $\nu=0.1$. We shall consider interior sensors with supports covering~$30\%$ of the volume of~$\Omega$ and boundary actuators covering~$30\%$ of the volume of the boundary~$\partial\Omega$. The configurations used in the simulations are shown in Fig.~\ref{fig:placements}.

We consider the following variation of~\eqref{J-Qz.v.u},
\begin{equation}\label{J-Qz.v.u-num}
J(y,u)\coloneqq\norm{Q(y_1+\clE\clB y_2)}{L^2(I;\bbR^s)}^2+q_v\norm{y_2}{L^2(I;\bbR^m)}^2+r\norm{u}{L^2(I;\bbR^m)}^2.
\end{equation}
The numerical costs are reported with the same normalization as the
abstract cost, without a common prefactor $1/2$. Accordingly, the cost
values obtained from the saved simulations with the former prefactor
$1/2$ are multiplied by two. This positive rescaling does not change the
exact minimizers or their associated trajectories; the reported state,
control, and estimation-error norms are unchanged.
The operator $Q$ is taken as in~\eqref{operBQ-same},
with a pair~$(q_v,r)\in\bbR_+^2$ of positive tuning parameters, which can play a relevant role on the performance of the input. In the simulations hereafter, we took~$(q_v,r)\coloneqq(10^{-8},10^{-5})$. For the actuators, we took the bump functions  
\[
\psi_i(\gamma)\coloneqq\indf_{\clI_i}\sin(\tfrac{\gamma-\clI_{i_1}}{\clI_{i_2}-\clI_{i_1}}\pi),\qquad \gamma\in\clI_i\coloneqq(\clI_{i_1},\clI_{i_2})\subset\partial\Omega,
\]
where the~$\clI_i$ are the intervals of the boundary (e.g., parametrized by arc-length) shown in Fig.~\ref{fig:placements}. Note that these actuators are in~$W^{1,2}(\partial\Omega)\subset W^{-\frac12,2}(\partial\Omega)$.
\subsection{Discretization}
We consider a triangulation of the spatial domain and use continuous $P_1$
finite elements (hat functions) to approximate our solution.  In time we use an IMEX scheme: the (implicit) Crank--Nicolson scheme for the diffusion operator (and linear reaction $ay$, for the linear model), and an explicit Adams--Bashforth extrapolation for the nonlinear term (and linear convection $b\cdot\nabla y$, for the linear model).
 The finite-horizon optimal control problems are solved with the alternating
Barzilai--Borwein method~\cite{BarzilaiBorwein88,AzmiKunisch22}, combined with a
nonmonotone line search; see, e.g.,~\cite{AzmiBern25}.

In all experiments, the time step is $\Delta t=0.005$. For each
actuator/sensor configuration, the mesh is generated with an initial
mesh-size parameter $h_{\max}=0.4$, with actuator endpoints and sensor
boundaries incorporated into the mesh. Two local refinements are applied
at each type of support interface, followed by one uniform refinement.
The optimizer uses absolute and initial-gradient-relative tolerances of
$10^{-8}$ and at most $500$ iterations. The nonmonotone line search uses
a memory of $30$ cost values, a sufficient-decrease parameter of $10^{-4}$,
and a backtracking factor of $1/2$.

\subsection{Nonlinear dynamics}\label{ssec:num-nonlinearRH}

We consider the Schl\"ogl equation in the unit square~$\Omega=(0,1)\times(0,1)\subset\bbR^2$, with boundary controls.
Here we consider the boundary-controlled system~\eqref{Schlogl-1},
 \begin{align}\label{Schlogl-num-bdry}
    &\tfrac{\p}{\p t} z=\nu\Delta z-z- \clP(t,z),\qquad
   (\bfn\cdot\nabla z)\rest{\p\Omega}=\sum_{i=1}^{m} v_i\psi_i,\qquad z(0)=z_0.
\end{align}    
with the Schl\"ogl nonlinearity in~\eqref{operSchlogl} and initial state $z_0=z_0(x)=2$. In particular, under Neumann boundary conditions, $z_0$ is a stable equilibrium and our asymptotic target~$y_\infty\coloneqq 0$ is an unstable equilibrium. Therefore, it is clear that a control input is necessary to stabilize the system to~$y_\infty$.

The receding horizon control~$v=v_{\rm rh}$ is computed through the variation of Algorithm~\ref{alg:RH} as follows in Algorithm~\ref{alg:RH-bdry}.   The reason we consider Algorithm~\ref{alg:RH-bdry} is motivated by eventual applications, where we are going to compute the input for the auxiliary extended model/system~\eqref{Schlogl-ext} and then apply it to the real system~\eqref{Schlogl-num-bdry}.

In the simulations, we took the prediction horizon and minimum sampling time as $T=1.5$ and $\tau=0.15$.

Recall the state space~$Y= V\times\bbR^m$, with~$V=W^{1,2}(\Omega)$. To underline that the norm of~$V$ depends on~$\nu$, in the following figures we denote
\begin{equation}
 |z|_{V_\nu}\coloneqq |z|_{V}
  =\left(
    |z|_{L^2(\Omega)}^2
    +\nu|\nabla z|_{L^2(\Omega)^2}^2
  \right)^{1/2},
  \label{eq:schloegl-state-norm}
\end{equation}
computed numerically as
$|z|_{V_\nu}^2\approx z_h^{\mathsf T}(\bfM+\nu \bfS)z_h$, where~$z_h$ is the finite-elements representation of~$z$ and~$\bfM$ and~$\bfS$ are the  mass and stiffness matrices.

For the squeezing-time selection, we use the following weighted norm
on~$Y$. Recall that~$v=y_2$ and that the physical state is reconstructed
as~$z=y_1+\clE\clB v$:
\begin{equation}
  |y|_{Y_{\nu,q_v}}\coloneqq\left(
    |y_1+\clE\clB v|_{V_\nu}^2
    +q_v|v|_{\bbR^m}^2
  \right)^{1/2}.
  \label{eq:schloegl-pair}
\end{equation}
Thus this norm measures the physical pair~$(z,v)$, with the same
positive weight~$q_v$ as above. Algorithm~\ref{alg:RH-bdry} uses this
norm when selecting the implemented prefix of each optimal prediction.

\begin{algorithm}[htbp]
\caption{Receding horizon boundary control,   RHC($\tau,T$)}\label{alg:RH-bdry}
\begin{algorithmic}[1]
\REQUIRE{lower bound $\tau$ for sampling time, prediction horizon $T\geq\tau$, a lifting operator~$E_{i,b}=\clE$ satisfying~\eqref{liftE}, initial data~$(z_0,v_0)$ with~$v_0\in\bbR^m$ and~$z_0-E_{i,b}\clB v_0\in V$}.
\ENSURE{a control input~$v_{\rm rh}\in W^{1,2}( \bbR_+ , \bbR^m)$; concatenation times $(t_k)_{k\in\bbN}$.}
\STATE Set~$\overline t_0\coloneqq 0$, $t_0=\overline t_0$, $k=0$, and~$\overline{y}_0\coloneqq(z_0-E_{i,b}\clB v_0,v_0)\in V\times\bbR^m$; initialize~$z_{\rm rh}(0)=z_0$.
\STATE\label{alg:BTinf} Set~$T_{\infty}=\infty$;
\WHILE{$\overline t_0<T_{\infty}$}
\STATE\label{alg:Bsolve} ``Solve'' the  open-loop problem~\ref{tag:opfin} for the extended system~\eqref{Schlogl-ext}, to find~$u^{\overline t_0,\overline{y}_0,*}_T  \in L^2(\clI_{\overline{t}_0,T},\bbR^m)$ and the corresponding solution~$ y=( y_1,y_2)$;
\STATE\label{alg:Bdelta} Set~$\delta:=\max\argmin\limits_{s\in[\tau,T]}\norm{y(\overline t_0+s)}{Y_{\nu,q_v}}^2$, using~\eqref{eq:schloegl-pair}.
\STATE Set~$t_{\rm new}\coloneqq\overline t_0+\delta$.
\STATE For~$t\in[\overline t_0,t_{\rm new}]$, set~$v_{\rm rh}(t)=y_2(t)$, including its continuous endpoint value.
\STATE Solve system~\eqref{Schlogl-num-bdry} from its current state~$z_{\rm rh}(\overline t_0)$, with~$v=v_{\rm rh}$, up to~$t_{\rm new}$.
\STATE Update~$\overline y_0\leftarrow\bigl(z_{\rm rh}(t_{\rm new})-E_{i,b}\clB v_{\rm rh}(t_{\rm new}),v_{\rm rh}(t_{\rm new})\bigr)$ and then~$\overline t_0\leftarrow t_{\rm new}$.
\STATE $k\gets k+1$;\quad $t_k\gets \overline t_0$.
\ENDWHILE
\end{algorithmic}
\end{algorithm}

\subsubsection{Actuator-sensor configurations: comparison on time interval \texorpdfstring{$[0,10]$}{[0,10]}}
We run the algorithm up to~$T_\infty=10$. Table~\ref{tab:schloegl-comparison} shows the norm of the state at final time~$t=10$ and the accumulated closed-loop cost. As we would expect, a larger number of actuators and sensors leads to a smaller norm at final time.
\begin{table}[!htbp]
  \centering
  \caption{Norm of final state and total cost.}
  \label{tab:schloegl-comparison}
  \begin{tabular}{rr||rr}
    \hline
    {$m$} & {$s$} & {$|z(10)|_{V_\nu}$} &  {$J_{[0,10]}(y^u,u)$}  \\
    \hline
    $4$  & $9$ & $1.7766$    &  $2.6752$  \\
    $12$ & $4$ & $0.43171$   &  $4.2756$  \\
    $12$ & $9$ & $0.0030709$ &  $2.2156$  \\
    \hline
  \end{tabular}
\end{table}

\begin{figure}[!htbp]
  \centering
  \includegraphics[width=\textwidth]{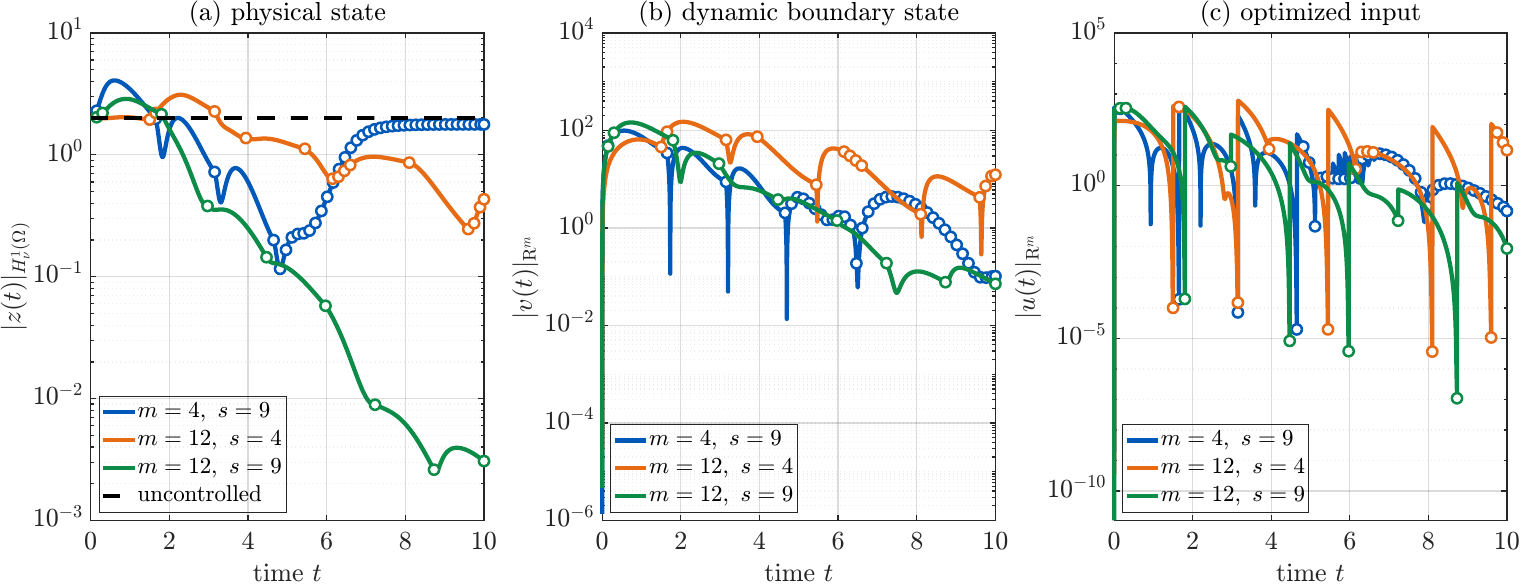}
  \caption{Evolution of the norms of the state~$|z(t)|_{V_\nu}$, the input~$|v(t)|_{\bbR^m}$, and the auxiliary/fictitious input~$|u(t)|_{\bbR^m}$.  Circle markers  are drawn at the RHC concatenation time instants found by Algorithm~\ref{alg:RH-bdry}.}
  \label{fig:schloegl-comparison}
\end{figure}
For the configurations tested with prediction horizon $T=1.5$, we observe stabilization with $m=12$ actuators and $s=9$ sensors, but not with $m=4$ actuators and $s=9$ sensors. With~$m=12$ actuators and~$s=4$ sensors, we see that the input has a poorer performance when compared with the case of~$s=9$ sensors. The stabilizing property of the input is also not clear enough; this suggests that~$4$ sensors may not be sufficient. Anyway, the results show that we obtain a better performance for larger numbers of actuators and sensors. Note that though this is expected, it is not just simple common sense, because the volume covered by (the supports of the) actuators and sensors is independent of~$(s,m)$.

\begin{figure}[!htbp]
  \centering
  \includegraphics[width=\textwidth]{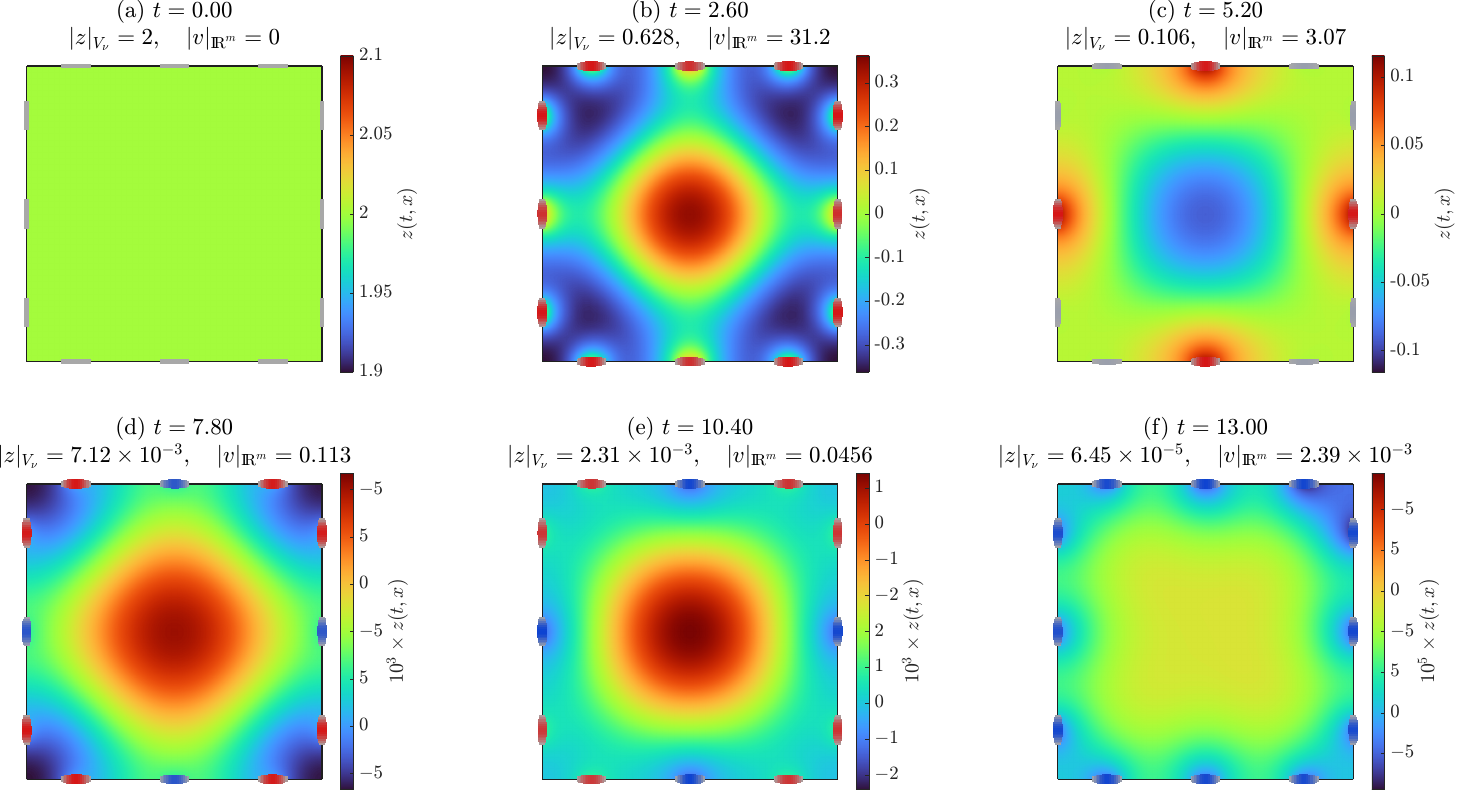}
  \caption{Physical state $z(t,x)$ for $s =9$ and $m =12$ at the shown time instants. The red
and blue colors at the boundary support segments of the actuators underline the positive and negative flux, respectively.}
  \label{fig:schloegl-snapshots}
\end{figure}

Figure~\ref{fig:schloegl-snapshots} shows snapshots of the state~$z(t,x)$ at
equally spaced time instants for $s=9$ and $m=12$.
Each snapshot uses an individual scaling of the state to preserve the visibility
of the spatial pattern over several orders of magnitude of decay. The positive
power of~$10$ displayed next to the colorbar indicates the corresponding
rescaling factor. As can be seen, the state decays to
$\norm{z(13)}{V_{\nu}} = 6.45\!\times\!10^{-5}$ at time $t=13$.

\subsection{Linear dynamics with dynamic output-feedback-inputs}\label{ssec:num-linearORHC}

Here we consider system~\eqref{Schlogl-num-bdry} with a linear reaction-convection term~$A_{\rm rc}z=A_{\rm rc}(x,t)z(x,t)$,
 \begin{align}\label{Schlogl-num-lin}
    &\tfrac{\p}{\p t} z=\nu\Delta z-z-A_{\rm rc}z,\qquad(\bfn\cdot\nabla z)\rest{\p\Omega}=\sum_{i=1}^{m} v_i\psi_{i},\qquad z(0)=z_0
\end{align}    
 with the receding horizon control~$v(t)=v_{\rm rh}(t)\in\bbR^{m}$, $m=12$, computed following a variant of Algorithm~\ref{alg:ORHC} based on  the extended system, as in Algorithm~\ref{alg:ORHC-bdry} below. 
 
 Here, we only consider the case of~$m=12$ actuators and~$s=9$ sensors.

Once the control input is chosen/computed, the observer for~\eqref{Schlogl-num-lin} is taken as
 \begin{align}\label{observer-num}
    &\tfrac{\p}{\p t} \widehat z=\nu\Delta \widehat z-\widehat z-A_{\rm rc}\widehat z+\fkI_s(Q\widehat z-w),\qquad(\bfn\cdot\nabla \widehat z)\rest{\p\Omega}=\sum_{i=1}^{m} v_i\psi_{i},\qquad \widehat z(0)=\widehat z_0
\end{align}  
where~$\widehat z_0$ is a guess/estimate we might have for~$z_0$.

Now, the state space for the extended system is~$Y= L^2(\Omega)\times\bbR^m$. Again we take a norm inspired by the cost functional,
\begin{equation}
  |y|_Y
  =\left(| y_1+\clE\clB y_2|_{L^2(\Omega)}^2
  +10^{-8}|y_2|_{\bbR^m}^2\right)^{1/2}.
  \label{eq:observer-pair}
\end{equation}

Here~$w=Qz$ is the output of the sensor measurements. The output-injection
operator~$\fkI_s\in\clL(\bbR^s,L^2(\Omega))$ is chosen to satisfy
(cf.~\cite{AzmiKunRod23,Rod21-aut})
\begin{equation}\label{out-inj-oper-num}
\fkI_sQh=-\lambda P_{\clW_s}h,\qquad h\in L^2(\Omega),
\end{equation}
where~$P_{\clW_s}\colon L^2(\Omega)\to\clW_s$ denotes the orthogonal
projection in~$L^2(\Omega)$ onto
$\clW_s\coloneqq\linspan\{\indf_{\omega_i^s}\mid 1\le i\le s\}$.
For the disjoint sensor supports used here, put
\[
 \mu_i\coloneqq|\omega_i^s|>0,\qquad
 \phi_i\coloneqq\frac{\indf_{\omega_i^s}}{\sqrt{\mu_i}},\qquad
 (Qh)_i=(h,\phi_i)_{L^2(\Omega)},\qquad 1\le i\le s.
\]
The normalized indicators~$\phi_i$ are orthonormal, so
$P_{\clW_s}h=\sum_{i=1}^s(Qh)_i\phi_i$. Consequently, the explicit
injection consistent with these normalized measurements is
\[
 \fkI_s q=-\lambda\sum_{i=1}^s q_i\phi_i
          =-\lambda\sum_{i=1}^s\frac{q_i}{\sqrt{\mu_i}}\indf_{\omega_i^s},
 \qquad q\in\bbR^s.
\]

\begin{remark}
Both the plant~$z$ and the observer~$\widehat z$ receive the same Neumann
boundary input~$\clB v$. Their difference
$e\coloneqq\widehat z-z$ therefore has homogeneous Neumann data and satisfies
\begin{align}\label{observer-num-error}
    &\tfrac{\p}{\p t} e=\nu\Delta e-e-A_{\rm rc}e-\lambda P_{\clW_s}e,\qquad(\bfn\cdot\nabla e)\rest{\p\Omega}=0,\qquad e(0)=\widehat z_0-z_0.
\end{align}
We know (e.g., by the arguments in~\cite{AzmiKunRod23,Rod21-aut}) that this system is exponentially stable, provided the number of sensors~$s$ (appropriately constructed, e.g., as in Fig.~\ref{fig:suppActGrid}) is large enough and provided the parameter~$\lambda\ge0$ is large enough as well.
 In particular, the assumption in eq.~\eqref{eq:exponent_output_inj} within Theorem~\ref{Thm:orhc-noq}  is satisfied. 
\end{remark}
\begin{algorithm}[htbp]
\caption{Output-based receding horizon boundary control, ORHC($\tau,T;\clZ,\fkI_S,E_{i,b},\widehat z_0$)}\label{alg:ORHC-bdry}
\begin{algorithmic}[1]
\REQUIRE{lower bound $\tau$ for sampling time; prediction horizon $T\ge \tau$; output operator $\clZ$; injection operator $\fkI_S$; a lifting operator~$E_{i,b}$ satisfying~\eqref{liftE}; initial guess/estimate $\widehat z_0\in L^2(\Omega)$ for the initial physical state~$z_0\in L^2(\Omega)$.}
\ENSURE{a control input $v_{\rm orh}\in L^2_{\rm ad}(\bbR_+,\bbR^m)$; observer/state estimate $\widehat z(\cdot)$; concatenation times $(t_k)_{k\in\bbN}$.}

\STATE Set $t_0\coloneqq0$, $\overline t_0\coloneqq0$, $\overline{y}_0\coloneqq(\widehat z_0,0)$ and~$\overline z_0=\widehat z_0$.
\STATE Set $T_\infty:=\infty$; $k:=0$.
\WHILE{$\overline t_0<T_\infty$}
  \STATE \label{alg:ORHC-bdry_Bsolve}``Solve'' the FTH problem~\ref{tag:opfin} for the extended system~\eqref{Schlogl-ext}, with $\clP(t,z)$ replaced by $A_{\rm rc}(t)z$, to obtain an optimal input
  $u_T^{\overline t_0,\overline{y}_0,*}\in L^2(\clI_{\overline t_0,T},\bbR^m)$ and the optimal predicted trajectory
  $y_T^{\overline t_0,\overline{y}_0,*}(t)\eqqcolon (\overline y_1(t),\overline y_2(t))\in Y$.
  \STATE \label{alg:ORHC-bdry_Bdelta} Set
  $\delta := \max\argmin_{s\in[\tau,T]} \norm{y_T^{\overline t_0,\overline{y}_0,*}(\overline t_0+s)}{Y}$.
  \STATE For all~$t\in [\overline t_0,\overline t_0+\delta]$, set $v_{\rm orh}(t):=\overline y_2(t)\in\bbR^m$;
 
  \STATE \label{alg:ORHC-bdry_runReal} Let the plant \eqref{Schlogl-num-lin} evolve on $[\overline t_0,\overline t_0+\delta)$ with input $v_{\rm orh}$ and record the output
  $w_{\rm orh}(t)\coloneqq Qz_{\rm orh}(t)\in\bbR^s$.
  
  \STATE Propagate the observer~\eqref{observer-num} on $[\overline t_0,\overline t_0+\delta)$ with~$(w,v,\widehat z(\overline t_0))=(w_{\rm orh}, v_{\rm orh},\widehat z_0)$.
 
  \STATE Update:
  $(\overline{t}_0,\widehat z_0, v_0)\leftarrow (\overline t_0+\delta,\widehat z(\overline t_0+\delta),v_{\rm orh}(\overline t_0+\delta))$.
   \STATE Update: $\overline{y}_0\leftarrow (\widehat z_0-E_{i,b}\clB v_0,v_0)$.
  \STATE $k\gets k+1$;\quad $t_k\gets \overline t_0$.
\ENDWHILE
\end{algorithmic}
\end{algorithm}

In the simulations we took $\lambda=30$. The initial state~$z_0$ is assumed to be unknown. To test the algorithm, we perform the simulations in the case that state was
\[
z_0=z(0,x)=1-2\cos(\pi x_1).
\]
The initial guess for the observer is taken as $\widehat z_0\coloneqq P_{\clW_S}z_0$. Note that~$\widehat z_0$ depends only on the available output~$w_0=Qz_0$ at initial time. Indeed, we have $P_{\clW_S}z_0=-\lambda^{-1}\fkI_Sw_0$.

\begin{remark}
    Note that in real-world applications, the physical system is~\eqref{Schlogl-num-lin}, while the observer~\eqref{observer-num} is an auxiliary system usually implemented in a machine. Therefore, the real vector to be input in the control system~\eqref{Schlogl-num-lin} is~$v$, whereas~$u$ is an auxiliary optimization variable, being the input of the auxiliary extended control system. Roughly speaking, $v$ is the real input, whereas~$u$ is an auxiliary/fictitious input.
\end{remark}

\begin{remark}
    In Algorithm~\ref{alg:ORHC-bdry}, Step~\ref{alg:ORHC-bdry_runReal}, the initial condition is $z_0=z(\overline t_0)\in z((t_k)_{k\in\bbN})$. These are (assumed) unknown states, but while performing our simulations, we have access and use these states in, and only in, Step~\ref{alg:ORHC-bdry_runReal} to simulate the real system. We underline that, within Algorithm~\ref{alg:ORHC-bdry}, these (unknown full) states are not used to compute the stabilizing input.
\end{remark}

We take $A_{\rm rc}y\coloneqq ay+b\cdot\nabla y$, with
the reaction and convection coefficients
\begin{align}
 a(t,x)\coloneqq -3.5-x_1-|\sin(6t+x_1)|,\qquad b(t,x)\coloneqq\begin{bmatrix}
     (x_1+x_2)\\ |\cos(6t)x_1x_2|
 \end{bmatrix}.
   \label{eq:linear-operator}
\end{align}

We ran Algorithm~\ref{alg:ORHC-bdry} up to time~$T_\infty=15$ and obtained the results in Table~\ref{tab:observer-result}. As we see, at the final time both the state and the state estimate error are relatively small.

\begin{table}[!htbp]
\centering
\caption{Initial and final state norms, estimation errors, and total accumulated cost.}
\label{tab:observer-result}
\begin{tabular}{@{}ccccc@{}}
\toprule
\multicolumn{2}{c}{State norm}
&
\multicolumn{2}{c}{Estimation error}
&
Total cost
\\
\cmidrule(lr){1-2}
\cmidrule(lr){3-4}
\cmidrule(l){5-5}
$|z(0)|_{L^2(\Omega)}$
&
$|z(15)|_{L^2(\Omega)}$
&
$|\widehat z(0)-z(0)|_{L^2(\Omega)}$
&
$|\widehat z(15)-z(15)|_{L^2(\Omega)}$
&
$J_{[0,15]}(y^u,u)$
\\
\midrule
$1.7319$
&
$4.5093\times10^{-4}$
&
$1.4086$
&
$1.8476\times10^{-6}$
&
$2509.0$
\\
\bottomrule
\end{tabular}
\end{table}

\begin{figure}[!htbp]
  \centering
  \includegraphics[width=\textwidth]{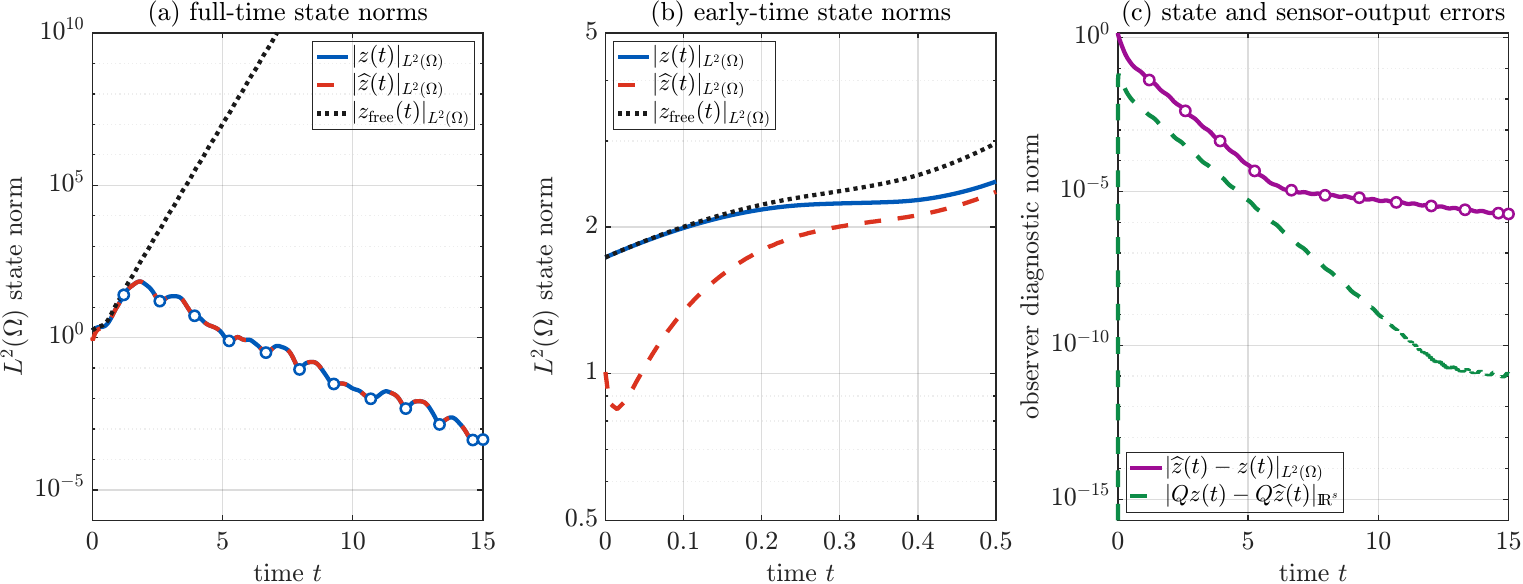}
  \caption{(a): Evolution of the norms of the uncontrolled state~$z_{\rm free}$, of the controlled state~$z$ and its estimate~$\widehat z$.  (b): Zoom of the evolution for small values of time. (c) evolution of the state estimate error~$\widehat z-z$ and of the output residual~$Q\widehat z-Qz$.}
  \label{fig:observer-performance}
\end{figure}

The free dynamics is unstable; quantitatively, we observed that the norm of the uncontrolled state reaches $\norm{z_{\rm free}(15)}{L^2(\Omega)}=2.09\times10^{21}$ at final time~$T_\infty=15$. This instability can also be observed in~Fig.~\ref{fig:observer-performance}. On the other hand, the results in the same figure show that the constructed input is able to stabilize the controlled state~$z$ to~$0$. Furthermore, the error of the state estimate also converges to~$0$, showing the stability of the closed-loop system.

\begin{figure}[!htbp]
  \centering
  \includegraphics[width=0.96\textwidth]{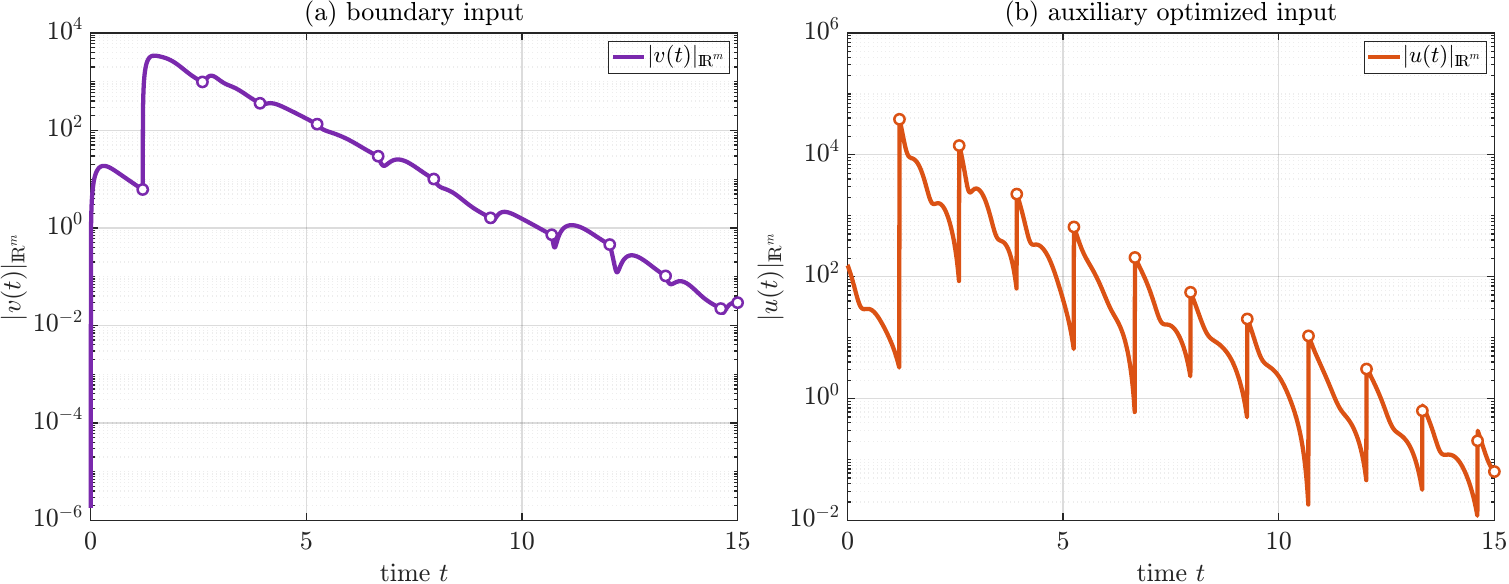}
  \caption{Evolution of the norms of the output-based boundary input~$v$ and of the optimized auxiliary/fictitious input~$u$. Circle markers are drawn at the ORHC concatenation time instants found by Algorithm~\ref{alg:ORHC-bdry}.}
  \label{fig:observer-boundary}
\end{figure}

Finally, Fig.~\ref{fig:observer-boundary} shows the evolution of the norms of
the boundary input~$v$ and the auxiliary input~$u$. Following the initial
transient, both norms decrease towards small values over the simulated
interval, consistently with the stabilizing behavior observed in
Fig.~\ref{fig:observer-performance}.

\subsection{Implementation of the boundary input}
The auxiliary input~$u$ is optimized over the prediction horizon~$T$
using the extended system~\eqref{Schlogl-ext}, with the corresponding
linear reaction--convection term in the linear example. The resulting
extended-state component~$v=y_2$ is then applied as the physical boundary
input to~\eqref{Schlogl-num-bdry} or~\eqref{Schlogl-num-lin} on
$[t_n,t_n+\delta_n)$, where~$\delta_n\in[\tau,T]$.

\section{Concluding Remarks}
\label{sec:conclusion}

In this paper, we developed an RHC (receding horizon control) framework for a class
of abstract semilinear evolution equations in continuous time. A central
feature of the analysis is that the stage cost is allowed to be only positive
semidefinite (i.e., not necessarily positive definite) with respect to the state. This situation arises naturally when
the state contribution to the cost is induced by finitely many sensors or
finite-dimensional observation operators. In such a setting, the cost does not
measure the full infinite-dimensional state, and stability has to be recovered
through suitable stabilizability and detectability-like mechanisms.

Under appropriate stabilizability and detectability-type assumptions, we
proved exponential stability of the RHC closed loop and derived a
suboptimality estimate with respect to the corresponding infinite-horizon
cost. In contrast to purely discrete-time formulations, the continuous-time
setting requires special care, since the relevant detectability properties are
finite-time properties and may depend on the length of the sampling interval.
This motivates the use of varying sampling times~$\delta_n$ that are selected online through the
synchronization mechanism of the proposed RHC scheme.

We also considered an ORHC (output-feedback RHC) framework in the linear case. While
the algorithmic structure follows our previous work, the stability analysis
presented here is different. By combining the RHC scheme with a
Luenberger-type state estimator and using an augmented Lyapunov functional, we
obtained a direct stability argument for the resulting observer-based closed
loop. This analysis highlights the role of the separation between the
estimator error dynamics and the RHC-stabilized state dynamics.

We applied the abstract framework to nonlinear Schl\"ogl equations with
finite-dimensional internal or dynamic boundary control, and to nonautonomous
linear parabolic equations in the output-feedback setting. In particular, the use of localized indicator-type functions for both
actuation and observation shows that the theory applies to systems in which
control and sensing are available only on small, distributed regions. The
numerical simulations support the theoretical findings and illustrate the
stabilizing behavior of the proposed full-state and ORHC
schemes.

Several directions remain open for future research. One natural question is
whether stability and suboptimality can be proven for a prescribed (fixed) sampling
time~$\delta_n=\tau$ under stronger uniform detectability assumptions, or by alternative
arguments; this question is important for improving robustness against small measurement/numerical errors. Another interesting direction is the extension of the
output-feedback analysis to nonlinear observer dynamics and more general
classes of semilinear parabolic equations; this is an important, but also a nontrivial long-standing open question in control design, due to the lack of the separation principle in the context of nonlinear dynamics. Finally, the combination of the
present framework with reduced-order or surrogate models, together with
a posteriori error certification would be important for reliable and
computationally efficient RHC implementations for large-scale PDE-constrained
control problems; here, the fact that the numbers of actuators and sensors are finite and the dimension of the range of~$Q$ is finite and, thus, independent of the reduced-order/surrogate models will plausibly be advantageous at the theoretical and numerical levels.
\subsection*{Acknowledgments}
The work of  S.R.  is funded by the FCT -- Portuguese Foundation for Science and Technology, I.P., under the scope of the projects UID/00297/2025~(\href{https://doi.org/10.54499/UID/00297/2025}{\nolinkurl{doi:10.54499/UID/00297/2025}});\; UID/PRR/00297/2025~(\href{https://doi.org/10.54499/UID/PRR/00297/2025}{\nolinkurl{doi:10.54499/UID/PRR/00297/2025}}); and UID/PRR2/00297/2025 (\href{https://doi.org/10.54499/UID/PRR2/00297/2025}{\nolinkurl{doi:10.54499/UID/PRR2/00297/2025}}) (Center for Mathematics and Applications -- NOVA Math).

After completing the initial draft, the authors used OpenAI Codex for
proofreading, consistency and citation checks, and assistance in reviewing
the mathematical arguments. The ideas, mathematical results and proofs,
algorithms, and numerical simulations are the authors' work, and the authors
take full responsibility for the entire content of the paper.

\appendix
\section{Well-posedness of the extended boundary-controlled Schl\"ogl system}
\label{app:schl-wellposedness}

This appendix provides the well-posedness and compactness properties used in
the verification of~\ref{asm:rh-existopt} and~\ref{asm:rh-estimYT} for Example~II.
The finite-horizon existence
of an optimal control is proved separately, in
Proposition~\ref{prop:schl-bdry-optimal-control}, because that argument is the
direct verification of~\ref{asm:rh-existopt}.
The proof uses the one-sided monotonicity of the Schl\"ogl nonlinearity. For local-in-time
existence, one may alternatively use the standard estimates in
\cite[Ass.~2.4]{Rod21-jnls}; see also
\cite[Ass.~3.4 and Sect.~5.2.1]{Rod20-eect}.

\begin{theorem}[Strong well-posedness]
\label{thm:schl-wellposedness}
Let $d\in\{1,2,3\}$, let $I=(\overline t_0,\overline t_0+T)$ be a nonempty
bounded interval, and suppose that the lifting
$\clE\clB:\bbR^m\to W^{2,2}(\Omega)$ is bounded. For every
$y_0\in Y=V\times\bbR^m$ and every $u\in L^2(I;\bbR^m)$, the extended
system~\eqref{Schlogl-ext} has a unique strong solution
\[
 y\in W(I;\rmD(A),H)\xhookrightarrow{}C(\overline I;Y).
\]
Moreover, for every $R>0$ there is a constant $C_{T,R}>0$ such that
\begin{equation}
\label{eq:schl-strong-bound}
 \norm{y}{C(\overline I;Y)}^2+
 \norm{y}{W(I;\rmD(A),H)}^2
 \leq C_{T,R}\left(\norm{y_0}{Y}^2+\norm{u}{L^2(I;\bbR^m)}^2\right)
\end{equation}
whenever
$\norm{y_0}{Y}^2+\norm{u}{L^2(I;\bbR^m)}^2\leq R^2$.
\end{theorem}

\begin{proof}\label{ssec:exist-sol-Schlogl}
Let $I=(\overline t_0,\overline t_0+T)$ and $u\in L^2(I;\bbR^m)$,
and write $L\coloneqq\clE\clB$. First solve
$\dot v+v=u$, $v(\overline t_0)=y_{20}$, exactly. Its solution belongs to
$W^{1,2}(I;\bbR^m)$ and satisfies
\begin{equation}\label{eq:galerkin-v-bound}
 \norm{v(t)}{\bbR^m}^2
 +\int_{\overline t_0}^{t}\norm{v(s)}{\bbR^m}^2\,\rmd s
 \le\norm{y_{20}}{\bbR^m}^2
       +\int_{\overline t_0}^{t}\norm{u(s)}{\bbR^m}^2\,\rmd s,
 \qquad t\in\overline I.
\end{equation}
Let $(e_j)_{j\ge1}$ be an orthonormal eigenbasis in~$L^2(\Omega)$ of
the homogeneous-Neumann operator~$A_0$. Put
$X_N=\linspan\{e_1,\ldots,e_N\}$ and let~$P_N$ denote the
$L^2(\Omega)$-orthogonal projection onto~$X_N$. These spectral projections
are contractions in~$L^2(\Omega)$ and~$V$, commute with~$A_0$ on its
domain, and satisfy~$P_Ny_{10}\to y_{10}$ in~$V$.
Define~$y^N=(y_1^N,v)$ by
\begin{equation}\label{eq:schl-galerkin-system}
 \dot y_1^N+A_0y_1^N+P_N\clP(y_1^N+Lv)=P_N f,
 \quad y_1^N(\overline t_0)=P_Ny_{10},
 \quad f\coloneqq L(v-u)-\overline A_0Lv.
\end{equation}
Here~$y_1^N(t)\in X_N$ and~$z^N\coloneqq y_1^N+Lv$.
The lifting assumption implies
$\norm{f(t)}{L^2(\Omega)}\le C(\norm{v(t)}{\bbR^m}
+\norm{u(t)}{\bbR^m})$, with~$C$ independent of~$N$.
The finite-dimensional system has a local solution by the Carath\'eodory
theorem. The estimates below extend it to all of~$I$.

In general~$z^N\notin X_N$, so we do not test the Galerkin equation
with~$z^N$ or directly transfer the physical-state energy identity to it.
Instead, we test~\eqref{eq:schl-galerkin-system} with~$2y_1^N$.
For~$\clP(r)=r^3-5r$, Young's inequality gives, for all~$r,a\in\bbR$,
\[
 -2\clP(r)(r-a)
 \le-r^4+15r^2+5a^2+\tfrac{27}{16}a^4.
\]
Applying this inequality with~$(r,a)=(z^N,Lv)$, and using
$\norm{Lv}{L^2(\Omega)}+\norm{Lv}{L^4(\Omega)}
\le C_1\norm{v}{\bbR^m}$, yields
\begin{align*}
 -2(\clP(z^N),y_1^N)_{L^2(\Omega)}
 &\le-\norm{z^N}{L^4(\Omega)}^4
       +C_2\norm{y_1^N}{L^2(\Omega)}^2
       +C_2\bigl(\norm{v}{\bbR^m}^2+\norm{v}{\bbR^m}^4\bigr),\\
 2(f,y_1^N)_{L^2(\Omega)}
 &\le\norm{y_1^N}{V}^2
       +C_3\bigl(\norm{v}{\bbR^m}^2+\norm{u}{\bbR^m}^2\bigr).
\end{align*}
Combining these inequalities with
$\tfrac{\rmd}{\rmd t}\norm{v}{\bbR^m}^2
+\norm{v}{\bbR^m}^2\le\norm{u}{\bbR^m}^2$, we obtain
\begin{equation}\label{eq:galerkin-basic-energy}
 \tfrac{\rmd}{\rmd t}\norm{y^N}{H}^2
 +\norm{y^N}{Y}^2+\norm{z^N}{L^4(\Omega)}^4
 \le C_4\norm{y^N}{H}^2
       +C_4\norm{v}{\bbR^m}^4+C_4\norm{u}{\bbR^m}^2,
\end{equation}
with~$C,C_i$, $1\le i\le4$,  all independent of
$\fkT\coloneqq(N,y_0,u)$ and of~$\overline t_0$.
For the dependence on the data, set
\begin{equation}
E_H^N\coloneqq\norm{y^N(\overline t_0)}{H}^2
                       +\norm{u}{L^2(I;\bbR^m)}^2,\quad
 E_Y^N\coloneqq\norm{y^N(\overline t_0)}{Y}^2
                       +\norm{u}{L^2(I;\bbR^m)}^2,\quad
 D_N\coloneqq E_H^N+(E_H^N)^2.\label{constENDN}
   \end{equation}
 
Estimate~\eqref{eq:galerkin-v-bound} implies
\[
 \int_I\norm{v(t)}{\bbR^m}^4\,\rmd t
 \le\norm{v}{L^\infty(I;\bbR^m)}^2
       \norm{v}{L^2(I;\bbR^m)}^2\le(E_H^N)^2.
\]
Gronwall's inequality and time integration of
\eqref{eq:galerkin-basic-energy} therefore bound
$\norm{y^N}{L^\infty(I;H)}^2+\norm{y^N}{L^2(I;Y)}^2$
by~$C(T)D_N$. These bounds prevent finite-time blow-up of the Galerkin
coefficients. Since~$z^N=y_1^N+Lv$ and~$L:\bbR^m\to V$ is bounded,
they also give
\begin{align}
 &\norm{y^N}{L^\infty(I;H)}^2
  +\norm{(z^N,v)}{L^\infty(I;H)}^2
  +\norm{y^N}{L^2(I;Y)}^2\notag\\
 &\qquad+\norm{z^N}{L^2(I;V)}^2
         +\norm{z^N}{L^4(I;L^4(\Omega))}^4
 \le C_5(T)D_N.
 \label{Galerkin-0}
\end{align}
Here~$C_5(T)$ is independent of~$\fkT$ and~$\overline t_0$.
The additional squared-data term in~$D_N$ accounts for the lifting in the
cubic nonlinearity; no uniform bound in~$N$ has been lost.

For the higher-order estimate, the test~$A_0y_1^N$ is admissible because
$A_0X_N\subseteq X_N$. Thus~$P_N$ drops out of its scalar product with
the nonlinear term. Together with the exact equation for~$v$, testing with
$2Ay^N$ gives
\begin{align}
\tfrac{\rmd}{\rmd t} \norm{y^N}{V\times\bbR^m}^2&=-2\norm{y^N}{\rmD(A)}^2-2(\clN(\cdot,y^N)-Bu,Ay^N)_{H}\notag\\
&\!\le -2\norm{y^N}{\rmD(A)}^2+2C_6\left(\norm{v}{\bbR^m}+\norm{u}{\bbR^m}\right)\norm{y^N}{\rmD(A)}-2(\clP(z^N),A_0y_1^N)_{L^2(\Omega)}. \label{Galerkin-1}
\end{align}
Next, we observe that
\begin{equation}\label{P-poly3-xi}
\clP(z^N)=(z^N)^3+\xi_2(z^N)^2+\xi_1z^N, 
\end{equation}
for constants~$\xi_1$ and~$\xi_2$ (in our example~$\xi_1=-5$ and~$\xi_2=0$).

Hence, denoting for simplicity $L^p\coloneqq L^p(\Omega)$ and using
the homogeneous Neumann condition on $y_1^N$, we find
\begin{align}
 -(\clP(z^N),A_0y_1^N)_{L^2}
 &=-\nu((3(z^N)^2+2\xi_2z^N+\xi_1)\nabla z^N,\nabla y_1^N)_{(L^2)^d}\notag\\
 &\quad-((z^N)^3+\xi_2(z^N)^2+\xi_1z^N,y_1^N)_{L^2}\notag\\
 &=-\nu((3(z^N)^2+2\xi_2z^N+\xi_1)\nabla z^N,\nabla z^N)_{(L^2)^d}\notag\\
 &\quad+\nu((3(z^N)^2+2\xi_2z^N+\xi_1)\nabla z^N,
                  \nabla\clE\clB v)_{(L^2)^d}\notag\\
 &\quad-((z^N)^3+\xi_2(z^N)^2+\xi_1z^N,z^N)_{L^2}\notag\\
 &\quad+((z^N)^3+\xi_2(z^N)^2+\xi_1z^N,\clE\clB v)_{L^2}.
 \label{test:PzAz1}
\end{align}
The lifting regularity gives
\[
 \norm{\nabla\clE\clB v}{(L^6)^d}
 +\norm{\clE\clB v}{L^6}\le C\norm{v}{\bbR^m},
\]
because $W^{2,2}(\Omega)\xhookrightarrow{}W^{1,6}(\Omega)$ for $d\le3$.
Both polynomials $-(3r^2+2\xi_2r+\xi_1)$ and
$-(r^2+\xi_2r+\xi_1)$ are bounded above.
Consequently, H\"older's inequality, with exponents $3,2,6$, gives
\begin{align}
 -2(\clP(z^N),A_0y_1^N)_{L^2}
 &\le C_5\norm{z^N}{V}^2
   +C_5\norm{3(z^N)^2+2\xi_2z^N+\xi_1}{L^3}
             \norm{z^N}{V}\norm{v}{\bbR^m}\notag\\
 &\quad+C_5\norm{(z^N)^2+\xi_2z^N+\xi_1}{L^3}
             \norm{z^N}{L^2}\norm{v}{\bbR^m}\notag\\
 &\le C_6\norm{z^N}{V}^2
   +C_7\left(\norm{z^N}{L^6}^2+\norm{z^N}{L^3}+1\right)
                   \norm{z^N}{V}\norm{v}{\bbR^m}\notag\\
 &\le C_7(1+\norm{z^N}{V}\norm{v}{\bbR^m})\norm{z^N}{V}^2
      +C_7\norm{z^N}{V}\norm{v}{\bbR^m}\notag\\
 &\le C_8(1+\norm{z^N}{V}^2+\norm{v}{\bbR^m}^2)\norm{z^N}{V}^2
      +C_8(\norm{z^N}{V}^2+\norm{v}{\bbR^m}^2)\notag\\
 &\le C_9(1+\norm{y^N}{V\times\bbR^m}^2)
                  \norm{y^N}{V\times\bbR^m}^2.
 \label{test:PzAz2}
\end{align}
which leads us to
\begin{align}
\tfrac{\rmd}{\rmd t} \norm{y^N}{V\times\bbR^m}^2&\le -2\norm{y^N}{\rmD(A)}^2+2C_6\left(\norm{v}{\bbR^m}+\norm{u}{\bbR^m}\right)\norm{y^N}{\rmD(A)}+C_9(1+\norm{y^N}{V\times\bbR^m}^2)\norm{y^N}{V\times\bbR^m}^2\notag\\
&\le -\norm{y^N}{\rmD(A)}^2+C_6^2\left(\norm{v}{\bbR^m}+\norm{u}{\bbR^m}\right)^2+C_9(1+\norm{y^N}{V\times\bbR^m}^2)\norm{y^N}{V\times\bbR^m}^2 \notag\\
&\le -\norm{y^N}{\rmD(A)}^2+C_{10}(1+\norm{y^N}{V\times\bbR^m}^2)\norm{y^N}{V\times\bbR^m}^2+C_{10}\norm{u}{\bbR^m}^2. \label{Galerkin-2}
\end{align}
The coefficient~$C_{10}(1+\norm{y^N}{Y}^2)$ in~\eqref{Galerkin-2}
has integral bounded by~$C(T)(1+D_N)$, by~\eqref{Galerkin-0}.
Gronwall's inequality, followed by time integration of~\eqref{Galerkin-2}
and enlargement of~$C_{11}(T)\ge1$, therefore gives
\begin{equation}\label{Galerkin-bddLinfVL2DA}
\norm{y^N}{L^\infty(I,V\times\bbR^m)}^2+\norm{y^N}{L^2(I,\rmD(A))}^2
\le C_{11}\rme^{C_{11}D_N}E_Y^N,
\end{equation}
with~$C_{11}=C_{11}(T)<\infty$ independent of~$\fkT$ and~$\overline t_0$.

Next, the polynomial estimate used for the time derivative is
\begin{align*}
 \norm{\clP(z^N)}{L^2(\Omega)}
 &\le C_{12}\left(\norm{z^N}{L^6(\Omega)}^3
        +\norm{z^N}{L^4(\Omega)}^2+\norm{z^N}{L^2(\Omega)}\right)\\
 &\le C_{13}\left(\norm{z^N}{V}^3+\norm{z^N}{V}^2+\norm{z^N}{V}\right)\\
 &\le C_{14}(1+\norm{y^N}{V\times\bbR^m}^2)
                     \norm{y^N}{V\times\bbR^m}.
\end{align*}
With~$E_H^N$, $E_Y^N$, and~$D_N$ defined above, in~\eqref{constENDN},  we put
$K_N\coloneqq C_{11}\rme^{C_{11}D_N}E_Y^N$.
The projection~$P_N$ is a contraction in~$L^2(\Omega)$, so it introduces
no dimension-dependent constant in the time-derivative estimate.
It follows from~\eqref{eq:schl-galerkin-system}, \eqref{Galerkin-0}, and
\eqref{Galerkin-bddLinfVL2DA} that
\begin{align*}
 \norm{\dot y^N}{L^2(I,H)}
 &\le C_{15}\bigl(\norm{y^N}{L^2(I,\rmD(A))}
           +\norm{v}{L^2(I,\bbR^m)}
           +\norm{u}{L^2(I,\bbR^m)}
           +\norm{\clP(z^N)}{L^2(I,L^2(\Omega))}\bigr)\\
 &\le C_{15}\bigl(\norm{y^N}{L^2(I,\rmD(A))}
           +\norm{v}{L^2(I,\bbR^m)}+\norm{u}{L^2(I,\bbR^m)}\bigr)\\
 &\quad+C_{15}(1+\norm{y^N}{L^\infty(I,V\times\bbR^m)}^2)
                      \norm{y^N}{L^2(I,V\times\bbR^m)}\\
 &\le C_{16}K_N^{1/2}+C_{16}(1+K_N)D_N^{1/2}\\
 &=C_{16}C_{11}^{1/2}\rme^{C_{11}D_N/2}(E_Y^N)^{1/2}
     +C_{16}(1+C_{11}\rme^{C_{11}D_N}E_Y^N)D_N^{1/2},
\end{align*}
where $C_{16}=C_{16}(T)$ is independent of $\fkT$. Combining the spatial
and time-derivative estimates gives
\begin{align}
 &\norm{y^N}{L^\infty(I,V\times\bbR^m)}^2
       +\norm{y^N}{W(I,\rmD(A),H)}^2 \le C_{17}\left(K_N+(1+K_N)^2D_N\right)\notag\\
 &=C_{17}C_{11}\rme^{C_{11}D_N}E_Y^N
       +C_{17}(1+C_{11}\rme^{C_{11}D_N}E_Y^N)^2D_N,
 \label{Galerkin-unifbound}
\end{align}
with $C_{17}=C_{17}(T)<\infty$ independent of $\fkT$, where
\[
 W(I,\rmD(A),H)
 \coloneqq\{h\in L^2(I,\rmD(A))\mid\dot h\in L^2(I,H)\}.
\]
In particular, if
$\norm{y_0}{Y}^2+\norm{u}{L^2(I;\bbR^m)}^2\le R^2$, then
$E_H^N\le E_Y^N\le R^2$ and
$D_N\le(1+R^2)E_Y^N$. Hence the right-hand side of
\eqref{Galerkin-unifbound} is at most~$C_{T,R}E_Y^N$, independently of~$N$.
After passage to the limit, the continuous embedding
$W(I,\rmD(A),H)\hookrightarrow C(\overline I;Y)$ yields precisely
\eqref{eq:schl-strong-bound}.

Since the right-hand side of~\eqref{Galerkin-unifbound} is bounded independently of~$N$, we may take a weakly converging subsequence (we still denote by~$y^N$),
\begin{align*}
y^N\xrightharpoonup[L^2(I,\rmD(A))]{}y^\infty,\qquad \dot y^N\xrightharpoonup[L^2(I,H)]{}\dot y^\infty.
\end{align*}
In particular, due to the compact embedding~$W(I,\rmD(A),H)\xhookrightarrow{\rm c}L^2(I,V\times\bbR^m)$ we obtain
\begin{align*}
y^N\xrightarrow[L^2(I,V\times\bbR^m)]{}y^\infty\quad\text{and}\quad z^N=y_1^N+\clE\clB v\xrightarrow[L^2(I,V)]{}z^\infty\coloneqq y_1^\infty+\clE\clB v.
\end{align*}
Note that~$\rmD(A)\xhookrightarrow{\rm d,c}V\times\bbR^m\xhookrightarrow{\rm d,c}H$, with compact and dense inclusions.

To show that~$y^\infty$ solves the dynamics in~\eqref{Schlogl-ext}, it is enough to show that each term converges weakly in~$L^2(I,L^2(\Omega))$. The convergence of the linear terms is straightforward. Thus, we focus on the nonlinear term, more precisely, on~$\clP(z^N)$. Recalling~\eqref{P-poly3-xi} and defining~$D^N\coloneqq z^N-z^\infty$,
\begin{align}
\clP(z^N)-\clP(z^\infty)&=(z^N)^3-(z^\infty)^3+\xi_2((z^N)^2-(z^\infty)^2)+\xi_1D^N\notag\\
&=D^N((z^N)^2+z^Nz^\infty+(z^\infty)^2)+\xi_2D^N(z^N+z^\infty)+\xi_1D^N. \label{Pz1-Pz2}
\end{align}
Given an arbitrary~$h\in L^2(I,L^2(\Omega))=L^2(I,L^2(\Omega))'$ we have that
\begin{align}
&\lim_{N\to\infty}\norm{(\clP(z^N)-\clP(z^\infty),h)_{L^2(I,L^2(\Omega))} }{}\notag\\
&\quad=\lim_{N\to\infty}\norm{(D^N((z^N)^2+z^Nz^\infty+(z^\infty)^2)+\xi_2D^N(z^N+z^\infty)+\xi_1D^N,h)_{L^2(I,L^2(\Omega))} }{}\label{decompPschl}\\
&\quad\le\lim_{N\to\infty}\norm{(D^N((z^N)^2+z^Nz^\infty+(z^\infty)^2),h)_{L^2(I,L^2(\Omega))} }{}\notag\\
&\quad\quad+\lim_{N\to\infty}\norm{(\xi_2D^N(z^N+z^\infty),h)_{L^2(I,L^2(\Omega))} }{}\notag\\
&\quad\quad+\lim_{N\to\infty}\norm{(\xi_1D^N,h)_{L^2(I,L^2(\Omega))}}{}.\notag
\end{align}
By the H\"older inequality, used in the form~$\norm{abc}{L^1(\Omega)}\le \norm{a}{L^6(\Omega)}\norm{b}{L^3(\Omega)}\norm{c}{L^2(\Omega)}$, we obtain
\begin{align*}
&\lim_{N\to\infty}\norm{(\clP(z^N)-\clP(z^\infty),h)_{L^2(I,L^2(\Omega))} }{}\\
&\qquad\le\lim_{N\to\infty}\norm{\norm{D^N}{L^6(\Omega)}\norm{(z^N)^2+z^Nz^\infty+(z^\infty)^2}{L^3(\Omega)}\norm{h}{L^2(\Omega)}}{L^1(I,\bbR)}\\
&\qquad\quad+\norm{\xi_2}{}\lim_{N\to\infty}\norm{\norm{D^N}{L^6(\Omega)}\norm{z^N+z^\infty}{L^3(\Omega)}\norm{h}{L^2(\Omega)}}{L^1(I,\bbR)}\\
&\qquad\quad+\norm{\xi_1}{}\lim_{N\to\infty}\norm{\norm{D^N}{L^2(\Omega)}\norm{h}{L^2(\Omega)}}{L^1(I,\bbR)}\\
&\qquad\le C_{18}\lim_{N\to\infty}\norm{\norm{D^N}{L^6(\Omega)}(\norm{z^N}{L^6(\Omega)}^2+\norm{z^\infty}{L^6(\Omega)}^2+1)\norm{h}{L^2(\Omega)}}{L^1(I,\bbR)}\\
&\qquad\le C_{19}\lim_{N\to\infty}\norm{\norm{D^N}{V}\norm{h}{L^2(\Omega)}}{L^1(I,\bbR)}\le C_{19}\lim_{N\to\infty}\norm{D^N}{L^2(I,V)}\norm{h}{L^2(I,L^2(\Omega))}=0.
\end{align*}
Therefore, $\clP(z^N)\xrightharpoonup[L^2(I,L^2(\Omega))]{}\clP(z^\infty)$.
To pass to the limit in the projected equation, note that for every
$h\in L^2(I;L^2(\Omega))$ we have~$P_Nh\to h$ in that space and
\[
 (P_N\clP(z^N),h)_{L^2(I;L^2(\Omega))}
 =(\clP(z^N),P_Nh)_{L^2(I;L^2(\Omega))}
\longrightarrow(\clP(z^\infty),h)_{L^2(I;L^2(\Omega))}.
\]
Also~$P_Nf\to f$ strongly.
We can conclude that~$y^\infty=(y_1^\infty,v)\in W(I,\rmD(A),H)$ is a solution of~\eqref{Schlogl-ext}.
The continuous trace map $W(I,\rmD(A),H)\to C(\overline I,Y)$ preserves the initial condition in the weak limit.
Since~$I=(\overline t_0,\overline t_0+T)$ with arbitrary~$\overline t_0$ and~$T>0$,
the solution is defined globally for~$t\in(\overline t_0,\infty)$
whenever $u\in L^2_{\rm loc}([\overline t_0,\infty),\bbR^m)$.

Finally, for uniqueness, let $y$ and $\overline y$ be two solutions for
the same $(y_0,u)\in Y\times L^2(I,\bbR^m)$, where $Y=V\times\bbR^m$ and
$y_0=(y_{10},y_{20})$. Their second components coincide, because both
satisfy $\dot v=-v+u$ with the same initial value. Explicitly,
\[
 y_2(t)=\overline y_2(t)\eqqcolon v(t)
 =\rme^{-(t-\overline t_0)}y_{20}
   +\int_{\overline t_0}^{t}\rme^{-(t-s)}u(s)\,\rmd s.
\]
Set $z=y_1+\clE\clB y_2$ and
$\overline z=\overline y_1+\clE\clB\overline y_2$.
Both $(z,v)$ and $(\overline z,v)$ solve~\eqref{Schlogl-ext0} with
$(z_0,v_0)=(y_{10}+\clE\clB y_{20},y_{20})$. Therefore the difference
$w=z-\overline z$ solves
\begin{align*}
 \tfrac{\p}{\p t}w&=\nu\Delta w-w-\clP(\cdot,z)+\clP(\cdot,\overline z),\\
 (\bfn\cdot\nabla w)\rest{\p\Omega}&=0,\qquad w(\overline t_0)=0.
\end{align*}
Testing with $2w$ gives
\begin{align*}
 \tfrac{\rmd}{\rmd t}\norm{w}{L^2(\Omega)}^2+2\norm{w}{V}^2
 &=-2(\clP(\cdot,z)-\clP(\cdot,\overline z),w)_{L^2(\Omega)}.
\end{align*}
Using the same decomposition as in~\eqref{decompPschl}, we obtain
\begin{align*}
 &-2(\clP(\cdot,z)-\clP(\cdot,\overline z),w)_{L^2(\Omega)}\\
 &\quad=-2(w(z^2+z\overline z+\overline z^2)
               +\xi_2w(z+\overline z)+\xi_1w,w)_{L^2(\Omega)}\\
 &\quad=-2\int_\Omega
       (z^2+z\overline z+\overline z^2+\xi_2(z+\overline z)+\xi_1)w^2\,\rmd x\\
 &\quad\le2\norm{z^2+z\overline z+\overline z^2
                 +\xi_2(z+\overline z)+\xi_1}{L^3(\Omega)}
                   \norm{w^2}{L^{3/2}(\Omega)}\\
 &\quad\le C_1\norm{w}{L^3(\Omega)}^2
          \le C_2\norm{w}{W^{1/2,2}(\Omega)}^2\\
 &\quad\le C_3\norm{w}{V}\norm{w}{L^2(\Omega)}
          \le\norm{w}{V}^2+C_4\norm{w}{L^2(\Omega)}^2.
\end{align*}
The constants are finite because both solutions are bounded in
$L^\infty(I,V)$. Hence
\[
 \tfrac{\rmd}{\rmd t}\norm{w}{L^2(\Omega)}^2+\norm{w}{V}^2
 \le C_4\norm{w}{L^2(\Omega)}^2.
\]
Gronwall's inequality and $w(\overline t_0)=0$ imply $w=0$.
Since $y_2=\overline y_2=v$, we also have $y_1=\overline y_1$, proving uniqueness.
\end{proof}

\bibliographystyle{plainurl}
\bibliography{RHC_OBS_Behzad}

\end{document}